\documentclass[reqno,11pt]{amsart}
\usepackage{amssymb}
\usepackage{verbatim}
\newif\ifpdf
\ifpdf
  \usepackage[pdftex]{graphicx}
  \usepackage[pdftex]{hyperref}
\else
  \usepackage{graphicx}
\fi
\numberwithin{equation}{section}       
\theoremstyle{plain}
\newtheorem{Thm}{Theorem}[section]
\newtheorem{Prop}[Thm]{Proposition}
\newtheorem{Lemma}[Thm]{Lemma}
\newtheorem{Cor}[Thm]{Corollary}
\newtheorem{Prop-def}[Thm]{Proposition-Definition}

\theoremstyle{definition}
\newtheorem{Example}[Thm]{Example}
\newtheorem{Def}[Thm]{Definition}
\newtheorem{Remark}[Thm]{Remark}

\newcommand{\C}{{\mathbf{C}}}

\newcommand{\N}{{\mathbf{N}}}
\newcommand{\Q}{{\mathbf{Q}}}
\newcommand{\R}{{\mathbf{R}}}
\newcommand{\Z}{{\mathbf{Z}}}

\newcommand{\fm}{{\mathfrak{m}}}

\newcommand{\bI}{{\bar{I}}}

\newcommand{\cA}{{\mathcal{A}}}

\newcommand{\cC}{{\mathcal{C}}}

\newcommand{\cH}{{\mathcal{H}}}

\newcommand{\cJ}{{\mathcal{J}}}
\newcommand{\cM}{{\mathcal{M}}}
\newcommand{\cN}{{\mathcal{N}}}
\newcommand{\cO}{{\mathcal{O}}}
\newcommand{\cP}{{\mathcal{P}}}

\newcommand{\cT}{{\mathcal{T}}}
\newcommand{\cV}{{\mathcal{V}}}

\newcommand{\hphi}{{\hat{\phi}}}

\renewcommand{\=}{ : = }
\renewcommand{\a}{\alpha}
\renewcommand{\b}{\beta}

\newcommand{\e}{\varepsilon}

\newcommand{\Spec}{\operatorname{Spec}}

\newcommand{\mass}{\operatorname{mass}}

\newcommand{\supp}{\operatorname{supp}}

\newcommand{\vol}{\operatorname{Vol}}
\renewcommand{\div}{\operatorname{div}}

\newcommand{\ie}{i.e.\ }

\newcommand{\qand}{{\quad\text{and}\quad}}

\newcommand{\vv}{{\vec{v}}}
\newcommand{\Lloc}{L^1_{\mathrm{loc}}}
\newcommand{\cVqm}{{\mathcal{V}_{\mathrm{qm}}}}
\newcommand{\cVdiv}{{\mathcal{V}_{\mathrm{div}}}}
\newcommand{\nuL}{{\nu^{\mathrm{L}}}}
\begin{document}
%
%
\setcounter{tocdepth}{1}

\title{Valuations and  multiplier ideals}
\date{\today}
\author{Charles Favre \and Mattias Jonsson}
\address{CNRS\\
         Institut de Math\'ematiques\\
         Equipe G\'eom\'etrie et Dynamique\\
         F-75251 Paris Cedex 05\\
         France}
\email{favre@math.jussieu.fr}
\address{Dept of Mathematics\\ 
  University of Michigan\\
  Ann Arbor, MI 48109-1109\\
  USA.}
\address{Dept of Mathematics\\ 
  Royal Institute of Technology\\
  SE-100 44 Stockholm\\
  Sweden.}
\email{mattiasj@umich.edu, mattiasj@kth.se}
\thanks{Second author partially supported by NSF Grant No DMS-0200614}
%
%
\begin{abstract}
  We present a new approach
  to the study of multiplier ideals in a local, two-dimensional 
  setting. Our method allows us to deal with
  ideals, graded systems of ideals and plurisubharmonic functions
  in a unified way. Among the applications are
  a formula for the complex integrability exponent of a 
  plurisubharmonic function in terms of  Kiselman numbers, 
  and a proof of the openness conjecture by Demailly and Koll\'ar.
  Our technique also yields new
  proofs of two recent results: one on the structure of
  the set of complex singularity exponents for holomorphic functions;
  the other by Lipman and Watanabe on the 
  realization of ideals as multiplier ideals.
\end{abstract}
%
%
\maketitle
%
%
\tableofcontents
%
%
%
%
%
%
\section*{Introduction}
This article is the third of a series of work on a new approach to
the study of singularities of various objects in a local, 
two-dimensional setting. Our focus in the present paper is on 
multiplier ideals and singularity exponents.

In the discussion below, we fix an
equicharacteristic zero, two-dimensional regular local ring $(R,\fm)$ 
with algebraically closed residue field. An important example is 
the ring $R=\cO_0$ of holomorphic germs at the origin in $\C^2$.

\smallskip
In~\cite{valtree}, we introduced the space $\cV$ consisting
of all $\R_+\cup\{+\infty\}$-valued valuations on $R$ 
centered at $\fm$, and normalized by $\nu(\fm)=1$. This space is
naturally a tree: it is a union of (uncountably many) real
segments patched together in such a way that $\cV$ remains homotopic
to a point. It is also an $\R$-tree in the classical sense for a
natural metric. We hence call $\cV$ the \emph{valuative tree}. 
It encodes in a natural way all
possible blowups of $R$ centered at $\fm$, and therefore gives a way
of measuring quite precisely singularities of different kinds of
objects. The points in $\cV$ that are not ends form the
subtree $\cVqm$ of \emph{quasimonomial} valuations. 
These valuations, which can alternatively be characterized 
as Abhyankar valuations of rank 1 or as pushforwards of
monomial valuations under birational morphisms, 
play a central role in our approach.

We used the valuative tree to study singularities of ideals 
$I\subset R$ in~\cite{valtree}, and 
of plurisubharmonic (psh) functions $u$ in~\cite{pshsing} 
(in which case $R=\mathcal{O}_0$). We showed that $I$ and $u$
both give rise to functions $g_I$ and $g_u$ on $\cVqm$, 
called the \emph{tree transforms} of $I$ and $u$.
These are defined as follows. For an ideal $I$ and $\nu\in\cVqm$,
we set $g_I(\nu):=\nu(I)=\min\{\nu(\phi)\ ;\ \phi\in R\}$. 
For a psh function, the definition of $g_u(\nu)=\nu(u)$ is more
involved but can be interpreted either as a generalized Lelong number 
in the sense of Demailly or as the pushforward of a Kiselman
number under a birational morphism. 

It turns out that the tree transforms $g_I$ and $g_u$ enjoy strong
concavity properties. We denote by $\cP$ the closure of the cone
generated by the functions on $\cVqm$ of the form
$c\,g_I$ over all $c>0$ and all ideals $I$ 
(or equivalently by all $g_u$'s for psh functions $u$).  An element in
$g\in\cP$ is called a \emph{tree potential}.  Tree potentials can
also be identified with positive measures on $\cV$. In the case of an
ideal, the measure is atomic and its decomposition into atoms is
equivalent to Zariski's factorization of integrally closed ideals.  In
the case of a psh function, the measure determines the most important
features of the singularity at the origin---a somewhat vague assertion
that the paper at hand serves to substantiate.  A general tree
potential may be loosely viewed as a formal analogue of a psh
function.

\smallskip
Multiplier ideals have emerged in recent years as a fundamental tool in
algebraic and analytic geometry. 
We refer to the book~\cite{laz-book}
for a detailed account on this subject, on its range
of applications, and for precise references.
See also~\cite{Dem} for a more analytic point of view.
Here we contend ourselves with recalling a few definitions 
adapted to our setting.

Fix a psh function $u$ defined near the origin in $\C^2$. 
The  multiplier ideal $\cJ(u) $ is the set of 
holomorphic germs $\phi\in\cO_0$ such that
$|\phi|^2\exp(-2u)$ is locally integrable. 
As $u$ may take the value
$-\infty$, $\cJ(u)$ is in general strictly included in
$\cO_0$. To worse behavior of $u$ at $0$, \ie for $u$
decreasing  faster to $-\infty$, corresponds deeper ideals $\cJ(u)$.
The multiplier ideal hence measures, in some sense, the
degree of singularity of $u$. 

To an ideal $I\subset R$ 
(or more generally a formal power $I^c$, $c>0$)
we can also associate a multiplier ideal $\cJ(I)$.
One way to do this uses resolution of singularities,
whereas Lipman~\cite{Lip94}
gives a more intrinsic definition.
It is also possible to associate an asymptotic multiplier ideal
to a graded system of ideals, \ie a sequence $(I_k)_{k=1}^\infty$ 
of ideals in $R$ such that $I_kI_l\subset I_{k+l}$.

\smallskip
Our first objective is to show that all these multiplier ideals
can be analyzed in a unified way using tree potentials. 
More precisely, to any tree potential $g$ on $\cVqm$
we associate a multiplier ideal 
$\cJ(g)\subset R$ as follows: 
an element $\psi\in R$ belongs to $\cJ(g)$ iff the function
$\chi_{\psi,g}(\nu)=g(\nu)/(\nu(\psi)+A(\nu))$ on $\cVqm$
is uniformly bounded by a constant $<1$. 
Here $A$ is a function on $\cVqm$ called \emph{thinness}.
It contains information on the relative canonical divisor 
of suitable birational models of $\Spec R$, 
as in the classical algebraic definition of multiplier ideals
(see Section~\ref{S6} for details).

We then show that the multiplier ideal of an ideal, a psh function,
or a graded system of ideals, coincides with that of the corresponding
tree transform. Let us be more precise.
If $I$ is an ideal and $c>0$, then we show that 
the multiplier ideal of $I^c$ coincides with that of the 
tree potential $c\,g_I$. 
Our proof essentially consists of translating
classical conditions for defining multiplier ideals into our language.
In the case of a graded system of ideals, the condition 
$I_kI_l\subset I_{k+l}$ implies a statement $g_{k+l}\le g_k+g_l$
on the level of tree potentials ($g_k$ is the tree 
potential of $I_k$). As we show, the sequence 
$k^{-1}g_k$ converges to a tree potential $g$ whose multiplier
ideal coincides with that of the graded system.

Similarly, we prove that if $u$ is psh, then the multiplier ideals 
$\cJ(g_u)$ and $\cJ(u)$ coincide. 
The proof is now more involved.
We first use Demailly's approximation technique
to reduce to the case of a psh function with logarithmic singularities.
The latter case is then reduced to the statement $\cJ(I^c)=\cJ(c\,g_I)$,
proved earlier, for a suitable chosen ideal $I$ and $c>0$.
Many of the arguments involved in the reductions
draw from~\cite{DK}, but our proof
also depends heavily on the fact that the tree transform $u\mapsto g_u$ 
behaves well under Demailly approximation.

\smallskip
Apart from providing a unifying framework, tree potentials
can be used as a powerful tool for studying many questions
regarding multiplier ideals. As support to this claim, we
give three applications of our approach.

The first two concern singularity exponents. To a general
tree potential $g$ is associated a number $c(g)$, called
the \emph{singularity exponent} or \emph{log-canonical
threshold} of $g$. It is defined by
$c(g)=\sup\{c>0\ ;\ \cJ(c\,g)=R\}$.
If $g=g_u$ is the tree transform of a psh function $u$, then
the singularity exponent is given by
$c(u)=\sup\{c>0\ ;\ \exp(-2cu)\in L^1_\mathrm{loc}\}$
and is a number measuring the ``strength'' of the singularity 
of $u$ at the origin. 
Various bounds were known for $c(u)$ in
terms of the Lelong number~\cite{skoda} and 
Kiselman numbers~\cite{kis2} of $u$. We sharpen these bounds, and
show that $c(u)$ can be in fact computed explicitly in terms 
of all Kiselman numbers of $u$.
We also prove that for $c=c(u)$, the function $\exp(-2cu)$ is
\emph{never} integrable. 
This provides an affirmative answer to the 
openness conjecture (in dimension 2) by Demailly and Koll\'ar
(see~\cite[Remark~5.3]{DK}). 
In fact, we establish the more precise estimate
$\vol\{u<\log r\}\gtrsim r^{2c(u)}$ as $r\to0$, and
we improve a recent result by Blel and Mimouni~\cite{BM,mim} 
by proving that
if $u$ has Lelong number one, then $\exp(-2u)$ fails to
be locally integrable at the origin iff $dd^cu$ is the sum 
of the current of integration on a smooth curve and a current
with zero Lelong number.

We then consider the set of complex singularity exponents
$\mathbf{c}=\{c(\log|\psi|)\}$ when $\psi$ ranges over all holomorphic
functions.  Shokurov~\cite{sho} used Mori's minimal program to show
that $\mathbf{c}$ satisfies the so-called \emph{ascending chain
condition}: any increasing sequence in $\mathbf{c}$ is eventually
stationary. Kuwata~\cite{Kuw99} subsequently improved this result by
characterizing explicitly all real numbers lying in $\mathbf{c}$. In
2000, Phong-Sturm in~\cite{phong-sturm} gave a completely analytic
proof of this result. We present an algebraic proof
in Section~\ref{sec-ACC}, independent of both approaches
above.  It is conjectured (see~\cite[Remark 3.5]{DK}) that the
ascending chain condition applies to set of complex singularity
exponents in any dimension (see~\cite{MP} for the most recent result
in this direction).  We hope that our result might lead to further
developments in higher dimension.

As a third application, we prove
that any integrally closed ideal in $R$ 
is the multiplier ideal of some formal power of an ideal.
This result was recently proved, independently, 
by Lipman and Watanabe~\cite{LipWat03}.

Finally we prove
quite generally that a tree
potential is completely characterized by the collection of multiplier
ideals $\{\cJ(tg)\}_{t\ge0}$. As a consequence, two psh functions
$u$ and $v$ have identical multiplier ideals $\cJ(tu)=\cJ(tv)$ 
for all $t\ge0$ iff they are \emph{equisingular} by which 
we mean that they have the same transforms $g_u=g_v$.
Equisingularity may be geometrically 
interpreted as follows:
for any composition $\pi$ of blowups,
the Lelong numbers of the pullbacks
$\pi^*u$ and $\pi^*v$ are the same at any point
on the exceptional divisor $\pi^{-1}(0)$, 
see~\cite[Proposition~6.2]{pshsing}.

\smallskip
Since the results on singularity exponents for psh functions
are arguably the most striking ones in the paper, 
we wish to briefly explain our approach, 
not using the language of valuations. 
Fix a psh function $u$ and assume, for simplicity,
that $dd^cu$ does not charge any curve.
To any irreducible (possibly singular) curve 
$D$ at the origin and any real number $t\ge 1$
we associate a family of punctured, conical neighborhood 
$\cA(r)$ of diameter $r$. Here $t$ determines
the ``thinness'' of the region. See Figure~\ref{F2} 
on p.\pageref{F2}. 
We have $\vol\cA(r)\sim r^{2A}$ for some $A>0$, and in
$\cA(r)$ we have $u\le\nu(u)\log r+O(1)$,
for some (maximal) real number $\nu(u)\ge0$.
If $\exp(-2cu)$ is locally integrable at the origin, 
it is integrable in $\cA(r)$, which easily implies 
$c\,\nu(u)<A$. 
The crucial fact is the existence of an \emph{optimal} 
region $\cA(r)$ detecting integrability: for this region 
we have $c(u)\nu(u)=A$, which yields nonintegrability of 
$\exp(-2c(u)u)$, \ie the openness conjecture. Moreover, the
curve $D$ associated to this region is \emph{smooth}, 
and then $\nu(u)$ is a Kiselman number of $u$. 
It is to establish the existence and main properties of 
the optimal region $\cA$ that we bring valuations into the picture.
In particular, the concavity properties of the tree transform
of $u$ play a key role.

\smallskip
Most results presented here are probably not particular to dimension
$2$. The main difficulty in extending our approach to higher
dimensions lies in understanding the analogue of the valuative tree.
We hope to tackle this problem in future research.

\smallskip
The present article relies in an essential way on the analysis
and formalism in our previous work: we recall in Section~\ref{S7} 
the main results of~\cite{valtree},~\cite{pshsing} that will be used.
The rest of the paper is then organized as follows: in 
Section~\ref{sec-def} we give the definition and main properties
of multiplier ideals of tree potentials. As we show in
Sections~\ref{sec-ideal} and~\ref{sec-psh},
this notion naturally generalizes multiplier
ideals of formal powers of ideals, graded system of ideals,
and psh functions. 
Sections~\ref{sec-arnold},~\ref{sec-ACC} 
and~\ref{sec-asmult} contain the applications mentioned above.
We study equisingularity in Section~\ref{sec-equi} and end the paper
with an appendix containing two proofs.
%
%
%
%
\section{Background}\label{S7}
In this section we give a brief review of the valuative tree and its
applications to the study of ideals and plurisubharmonic (psh)
functions.  For details, we refer to~\cite{valtree} for Sections~\ref{S13}
to~\ref{S5},  and to~\cite{pshsing} for Section~\ref{S20}.
%
%
\subsection{Conventions}\label{S13}
In general, $(R,\fm)$ denotes an
equi\-characteristic zero, two-dimensional 
regular local ring with algebraically closed 
residue field $k$. We will refer to this as the \emph{general case}.

Whenever we talk about psh functions, we will always 
be in the \emph{analytic case}, meaning that 
$R=\cO_0$ is the ring of holomorphic germs at the origin in $\C^2$.
Then $k=\C$ and $\fm$ is the maximal ideal of germs 
vanishing at the origin. 

In general we write $(\hat{R},\hat{\fm})$ for the completion of $R$.
It is the ring of formal power series in two variables with 
coefficients in $k$.
%
%
\subsection{The valuative tree}\label{S12}
Our starting point is the approach to valuations
worked out in~\cite{valtree}.
%
%
\subsubsection{Valuations} 
(\cite[Section~1.2]{valtree}.)
We consider the space $\cV$ of centered, normalized valuations on $R$,
\ie the set of functions $\nu:R\to[0,\infty]$ satisfying:
\begin{enumerate}
\item[(i)]\label{p1} 
  $\nu(\psi\psi')=\nu(\psi)+\nu(\psi')$ for all $\psi,\psi'$; 
\item[(ii)]\label{p2}
  $\nu(\psi+\psi')\ge\min\{\nu(\psi),\nu(\psi')\}$ for all
  $\psi,\psi'$; 
\item[(iii)]\label{p3} 
  $\nu(0)=\infty$, 
  $\nu|_{\C^*}=0$, 
  $\nu(\fm):=\min\{\nu(\psi)\ ;\ \psi\in\fm\}=1$.
\end{enumerate} 
Then $\cV$ is equipped with a natural \emph{partial ordering}: 
$\nu\le\mu$ iff $\nu(\psi)\le\mu(\psi)$ for all $\psi\in\fm$. 
The \emph{multiplicity valuation} $\nu_\fm$ defined by
$\nu_\fm(\psi)=m(\psi)=\max\{k\ ;\ \psi\in\fm^k\}$ 
is the unique minimal element of $\cV$. 
(See~\cite[Section~1.5.1]{valtree}.)

Any valuation on $R$ extends uniquely to a valuation in its
completion $\hat{R}$, hence the valuation spaces attached to
$R$ and $\hat{R}$ are isomorphic. (See~\cite[Theorem~3.1]{spiv}.)
%
%
\subsubsection{Curve valuations}\label{curve-val}
(\cite[Section~1.5.5]{valtree}.)
Some natural maximal elements are the \emph{curve valuations} defined
as follows. To each irreducible (possibly formal) curve $C$
we associate $\nu_C\in\cV$ defined by
$\nu_C(\psi)=C\cdot\{\psi=0\}/m(C)$, where ``$\cdot$'' denotes
intersection multiplicity and $m$ multiplicity.
If $C$ is defined by $\phi\in\hat\fm$, then we also write 
$\nu_C=\nu_\phi$. Note that $\nu_\phi(\phi)=\infty$.

The set $\cC$ of local irreducible curves carries a natural
(ultra)metric in which $\cC$ has diameter 1.
It is given by $d_\cC(C,D)=m(C)m(D)/C\cdot D$.
(See~\cite[Lemma~3.56]{valtree}.)
%
%
\begin{figure}[ht]
  \includegraphics[width=0.85\textwidth]{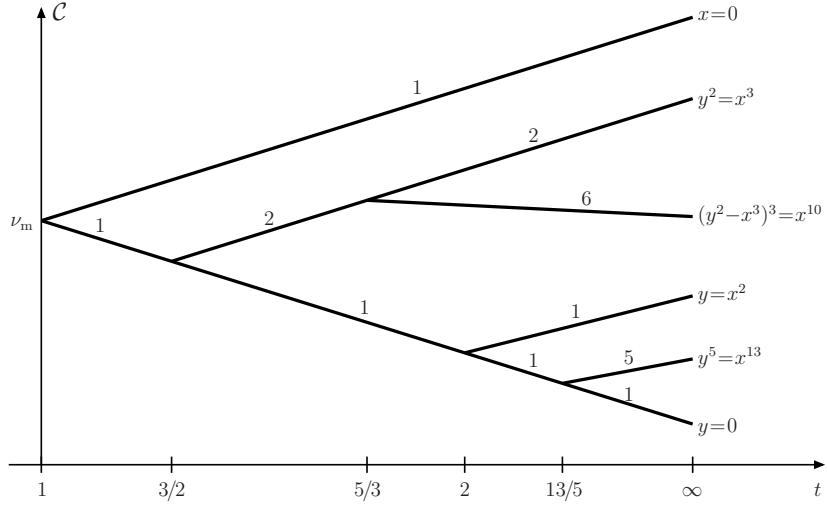}
  \caption{The valuative tree. The segments consist of valuations
    of the form $\nu_{\phi,t}$, where $\phi=x,y^2-x^3,\dots,y$
    and the skewness parameter $t$ ranges from $1$ to $\infty$.
    Skewness $t=1$ gives the multiplicity valuation $\nu_\fm$ and
    skewness $t=\infty$ the curve valuation $\nu_\phi$.
    The integer label above a segment indicates multiplicity
    (Section~\ref{S4}).}\label{F1}
\end{figure}
%
%
\subsubsection{Quasimonomial valuations}\label{sec-qm}
(\cite[Section~1.5.4]{valtree}.)
Arguably the most important valuations in $\cV$ are the
\emph{quasimonomial} ones.\footnote{A quasimonomial valuation can be
made monomial (\ie completely determined by its values on a pair of
local coordinates $(x,y)$) by a birational morphism. Quasimonomial
valuations are also known as Abhyankar valuations of rank 1,
see~\cite{ELS}.}  
They are of the form $\nu_{C,t}$, where $C\in\cC$ and $t\in[1,\infty)$, 
and satisfy $\nu_{C,t}(\psi)=\min\{\nu_D(\psi)\ ;\ d_\cC(C,D)\le
t^{-1}\}$.  We have $\nu_{C,s}=\nu_{D,t}$ 
iff $s=t\ge d_\cC(C,D)^{-1}$.
(See~\cite[Proposition~3.55]{valtree}.)
Thus $\cVqm$, the set of all
quasimonomial valuations, is a naturally a quotient of
$\cC\times[1,\infty)$, and has a natural tree structure: if
$\nu,\nu'\in\cVqm$ and $\nu<\nu'$, then the \emph{segment}
$[\nu,\nu']=\{\mu\in\cVqm\ ;\ \nu\le\mu\le\nu'\}$ is isomorphic to a
compact real interval.
See~\cite[Theorem~3.57]{valtree} and Figure~\ref{F1}.  
We set $\nu_{\phi,t}:=\nu_{C,t}$ when
$C=\{\phi=0\}$.

Quasimonomial valuations are of two types: \emph{divisorial}
and \emph{irrational}, depending on whether the parameter $t$
is rational or irrational.\footnote{A quasimonomial valuation $\nu$ 
  is irrational iff $\nu(R)\not\subset\Q$, hence the name.}
We write $\cVdiv$ for the set of divisorial valuations.

The full space $\cV$ is the completion of $\cVqm$ in the sense
that every element in $\cV$ is the limit of an increasing sequence
in $\cVqm$. It is hence also naturally a tree, called the 
\emph{valuative tree}. The ends of $\cV$ are exactly the elements
of $\cV\setminus\cVqm$ and are either curve valuations or
\emph{infinitely singular} valuations.\footnote{The latter are
represented by infinite Puiseux series whose exponents are
rational numbers with unbounded denominators.}
%
%
\subsubsection{Skewness and intersection multiplicity}\label{sec-skewness}
(\cite[Section~3.3]{valtree}.)
An important invariant of a valuation is its \emph{skewness} $\a$
defined by $\a(\nu)=\sup\{\nu(\phi)/m(\phi)\ ;\ \phi\in\fm\}$.
Skewness naturally \emph{parameterizes} the trees $\cVqm$ and $\cV$ in
the sense that $\a:\cVqm\to[1,\infty)$ is strictly increasing and
restricts to a bijection onto its image on any segment; indeed
$\a(\nu_{\phi,t})=t$ for any $\nu_{\phi,t}\in\cVqm$. 
Thus divisorial (irrational) valuations have rational
(irrational) skewness. Curve valuations have infinite skewness whereas
the skewness of an infinitely singular valuations may or may not be
finite.

The tree structure on $\cV$ implies that any collection $(\nu_i)_{i\in I}$
of valuations in $\cV$ admits an 
\emph{infimum} $\wedge_i\nu_i$, see~\cite[Corollary~3.15]{valtree}.
Together with skewness, this allows us to define an 
\emph{intersection product} on $\cV$:
we set $\nu\cdot\mu:=\a(\nu\wedge\mu)\in[1,\infty]$. 
This is a normalized extension of the intersection product on $\cC$ 
as $C\cdot D=(\nu_C\cdot\nu_D)m(C)m(D)$.
If $\nu\in\cV$ and $\phi\in\fm$ is irreducible, then 
$\nu(\phi)=m(\phi)\,(\nu\cdot\nu_\phi)$.
Moreover, if $\nu(\phi)$ is irrational, then 
$\nu=\nu_{\phi,t}$ with $t=\nu(\phi)/m(\phi)$.
%
%
\subsubsection{Tangent space and weak topology}
(\cite[Section~3.1-3.2]{valtree}.)
Let $\mu$ be a valuation in $\cV$. Declare 
$\nu,\nu'\in\cV\setminus\{\mu\}$ to be equivalent if
the segments $]\mu,\nu]$ and $]\mu,\nu']$ intersect.
An equivalence class is called a \emph{tangent vector} at $\mu$
and the set of tangent vectors at $\mu$, the \emph{tangent space},
denoted by $T\mu$. If $\vv$ is a tangent vector, we denote by
$U(\vv)$ the set of points in $\cV$ defining the equivalence class
$\vv$. The points in $U(\vv)$ are said to \emph{represent} $\vv$.

A point $\mu$ in the tree $\cV$ is an \emph{end}, 
a \emph{regular point},
or a \emph{branch point} when $T\mu$ contains one,
two, or three or more points, respectively.  In terms of valuations:
the ends of $\cV$ are curve and infinitely singular valuations;
the regular points irrational valuations; and the branch points
divisorial valuations, at which the tangent space is in bijection with
the complex projective line $\mathbf{P}^1$ and hence
uncountable. 
See~\cite[Theorem~3.20]{valtree}.

We endow $\cV$ with the \emph{weak topology},
generated by the sets $U(\vv)$ over all tangent vectors $\vv$;
this turns $\cV$ into a compact (Hausdorff) space.
If $\nu_k\to\nu$, 
then $\nu_k\wedge\mu\to\nu_k\wedge\mu$ for all $\mu\in\cV$.
The weak topology on $\cV$ is
characterized by $\nu_k\to\nu$ iff $\nu_k(\phi)\to\nu(\phi)$
for all $\phi\in R$, see~\cite[Theorem~5.1]{valtree}.
%
%
\subsubsection{Multiplicities}\label{S4}
(\cite[Section~3.4]{valtree}.)
By setting $m(\nu)=\min\{m(C)\ ;\ C\in\cC,\ \nu_C\ge\nu\}$ we
extend the notion of \emph{multiplicity} from $\cC$ to $\cVqm$.  
Clearly $m:\cVqm\to\N$ is increasing and hence extends to all of $\cV$. 
In fact $m(\nu)$ divides $m(\mu)$ whenever $\nu\le\mu$.  The infinitely
singular valuations are characterized as having infinite multiplicity.

As $m$ is increasing and integer valued, it is piecewise 
constant on any segment $[\nu_\fm,\nu_\phi]$, where $\phi\in\cC$.
This implies that $m(\vv)$ is naturally defined for any tangent
vector $\vv$. If $\nu$ is nondivisorial, then $m(\vv)=m(\nu)$ 
for any $\vv\in T\nu$. 

If $\nu$ is divisorial, then the situation is more complicated.
Suffice it to say that there exists an integer $b(\nu)$,
divisible by $m(\nu)$, such that $m(\vv)=b(\nu)$ for all
but at most two tangent vectors $\vv$ at $\nu$. We call
$b(\nu)$ the \emph{generic multiplicity} of $\nu$.
%
%
\subsubsection{Approximating sequences}\label{S3}
(\cite[Section~3.5]{valtree}.)
Consider a quasimonomial valuation $\nu\in\cVqm$. The
multiplicity $m$ is integer-valued and piecewise constant on the
segment $[\nu_\fm,\nu]$, hence has a finite number $g$
(possibly zero) of jumps.  Thus there are divisorial valuations
$\nu_i$, $0\le i\le g$ and integers $m_i$, such that
\begin{equation}\label{e705}
  \nu_\fm=\nu_0<\nu_1<\dots<\nu_g<\nu_{g+1}=\nu
\end{equation}
and $m(\mu)=m_i$ for $\mu\in\,]\nu_i,\nu_{i+1}]$, $0\le i\le g$.
We call the sequence $(\nu_i)_{i=1}^g$ the 
\emph{approximating sequence} associated to $\nu$.
It plays a prominent role in~\cite{spiv}. 

The concept of approximating sequences extends naturally to valuations
that are not quasimonomial: for curve valuations the sequences
are still finite, for infinitely singular valuations they are infinite.
%
%
\subsubsection{Thinness}\label{thinness}
(\cite[Section~3.6]{valtree}.)
Skewness $\a$ is a parameterization of $\cV$ that does not ``see'' 
multiplicities. Another parameterization, of crucial importance,
is \emph{thinness} $A$, defined as follows. 
If $\nu\in\cVqm$ then
\begin{equation}\label{e40}
  A(\nu)=2+\int_{\nu_\fm}^\nu m(\mu)\,d\alpha(\mu).
\end{equation}
In terms of~\eqref{e705} we have $A(\nu)=2+\sum_0^gm_i(\a_{i+1}-\a_i)$
with $\a_i=\a(\nu_i)$. Note that $A(\nu)\le 1+ m(\nu)\a(\nu)$.  Just
like skewness, we may define $A(\nu)$ also for $\nu\notin\cVqm$.

The inequality $A(\nu)\ge 1+\a(\nu)$ always holds, with equality iff
$m(\nu)=1$. Moreover, $A(\nu)<m(\nu)\a(\nu)$ as soon as $m(\nu)>1$,
and $A(\nu)-m(\nu)\a(\nu)\to-\infty$ as $\nu$ increases towards 
an infinitely singular valuation. 
%
%
\subsubsection{Geometric interpretation of divisorial valuations}\label{S6}
(\cite[Chapter~6]{valtree}.)
Every divisorial valuation $\nu$ arises as follows: there exists a
finite composition of point blowups $\pi:X\to\Spec R$ (\ie
$\pi:X\to(\C^2,0)$ in the analytic case) and an 
\emph{exceptional component} $E$ (\ie an irreducible component
of $\pi^{-1}(0)$) such that $\nu=b^{-1}\pi_*\div_E$, where
$b=b(\nu)$ is the generic multiplicity at $\nu$ and $\div_E$ denotes
the order of vanishing along $E$. 
The generic multiplicity of $\nu$
is then equal to the multiplicity of any \emph{curvette} of $\nu$,
\ie the image by $\pi$ of any smooth curve intersecting $E$ 
transversely at a smooth point. 
In fact, every generic tangent vector at $\nu$  is represented by
the curve valuation associated 
to a curvette~\cite[Section~6.6.1]{valtree}.
We also have $A(\nu)=a/b$, where $a-1$ is equal to
the order of vanishing along $E$ of the Jacobian determinant 
of $\pi$~\cite[Theorem~6.22]{valtree}.
%
%
\subsubsection{Borel measures and tree potentials}\label{sec-measure}
(\cite[Chapter~7]{valtree} or~\cite[Section~5]{pshsing}.)
Let $\cM$ be the space of (weak, regular) positive Borel measures 
on $\cV$,\footnote{In~\cite{valtree}, $\cM$ is denoted by $\cM^+$.}
endowed with the topology of vague convergence. 
We will identify $\nu\in\cV$ with its Dirac mass
$\delta_\nu\in\cM$, thus viewing $\cV$ as a (compact) subset of $\cM$.

Any $\rho\in\cM$ determines a real-valued function $g_\rho$ on $\cVqm$
defined by $g_\rho(\nu)=\int_\cV\mu\cdot\nu\,d\rho(\mu)$.  A function
of the form $g_\rho$ is called a \emph{tree potential}.
The space $\cP$ of tree potentials
is a closed convex cone in $\cV_{\mathrm{qm}}^{\R}$
(under pointwise convergence) and the map $\rho\mapsto g_\rho$ is a
homeomorphism of $\cM$ onto $\cP$ whose inverse naturally defines a
\emph{Laplace operator} 
$\Delta:\cP\to\cM$.\footnote{In~\cite{valtree}, $\cP$
  is denoted by $\cP^+$, and tree potentials are called 
  \emph{positive tree potentials}.} 
The mass of $\Delta g$ is $g(\nu_\fm)$. 
See~\cite[Theorem~7.64]{valtree}.

If $g\in\cP$, then $g$ is increasing and concave on any segment
$[\nu_\fm,\nu_0]$ (parameterized by skewness) in $\cV$. 
If $\nu\in\,]\nu_\fm,\nu_0[$, then the 
left derivative of $g$ at $\nu$
with respect to skewness is equal to 
$\rho\{\mu\ge\nu\}$,
whereas the right derivative is given by
$\rho (U(\vv))$, where $\rho=\Delta g$ and 
$\vv$ denotes the tangent vector at $\nu$ 
represented by $\nu_0$.
We always have 
$\rho\,\{\mu\ge\nu\}\a(\nu)\le g(\nu)\le g(\nu_\fm)\a(\nu)$ with 
equality (in either inequality) iff $\rho=\Delta g$ 
is supported on $\{\mu\ge\nu\}$.

A \emph{subtree} of $\cV$ is a subset $\cT$ such that $\nu\in\cT$
and $\mu\le\nu$ implies $\mu\in\cT$. A subtree is \emph{finite}
if it has finitely many ends.
If $g\in\cP$ and $\cT$ is
a subtree of $\cV$, then $g_\cT$ denotes the infimum of all tree
potentials coinciding with $g$ on $\cT$. 

The \emph{support} of
$g\in\cP$ is the smallest subtree $\cT$ for which $g=g_\cT$.
Alternatively, it is the smallest subtree containing $\supp\Delta g$.
The support of any tree potential is included in the closure of  
a countable union of finite trees.
%
%
\subsection{Ideals}\label{S5}
(\cite[Section~8.1]{valtree}.)
An ideal $I$ in $R$ is \emph{primary} if $\fm^n\subset I$ for
some $n>0$. The \emph{integral closure} of $I$ is the ideal $\bI$ of 
$\phi\in R$ such that $\phi^n+a_1\phi^{n-1}+\dots+a_n=0$ for some 
$n\ge1$ and $a_i\in I^i$. If $I=\bI$, then
$I$ is \emph{integrally closed}.

\medskip
Any ideal $I$ in $R$ has an associated \emph{tree transform} $g_I$,
defined by $g_I(\nu)=\nu(I)=\min\{\nu(\phi)\ ;\ \phi\in R\}$. This
function $g_I$ belongs to $\cP$, and thus defines a \emph{tree
measure} $\rho_I=\Delta g_I\in\cM$.  The latter measure has mass
$m(I):=\nu_\fm(I)$.

If $I=\phi R$ is principal, then we write $g_\phi=g_I$ 
and $\rho_\phi=\rho_I$: the latter measure is given by
$\rho_\phi=\sum_in_im_i\nu_i$, where 
$\phi=\prod\phi_i^{n_i}$ is the factorization of $\phi$ into
irreducible factors, $m_i=m(\phi_i)$ and $\nu_i$ the curve
valuation associated to $\phi_i$. 

If $I$ is primary, then
$\rho_I=\sum_in_ib_i\nu_i$, where $n_i\in\N$ and
$\nu_i$ are divisorial valuations with generic multiplicity $b_i$.
The valuations $\nu_i$ are exactly the \emph{Rees valuations} of $I$.

A general ideal $I$ is the product of a principal ideal and a 
primary ideal, hence has a tree measure of the form
$\rho_I=\sum_in_ib_i\nu_i+\sum_in_im_i\nu_i$. 

To any measure $\rho$ is associated an ideal $I_\rho=\{\phi\in R\ ;\
g_\phi\ge g_\rho\}$.  When $\rho$ is of the previous form, then
$\rho_{I_\rho}=\rho$, whereas $I_{\rho_I}$ is the integral closure of
$I$.  The decomposition of $\rho_I$ above gives an interpretation of
Zariski's factorization theorem: if $I$ is integrally closed, then
$I=\prod_iI_{\nu_i}^{n_i}$.  If furthermore $I$ is primary, \ie all
the $\nu_i$ are divisorial, then $I=\bigcap_i\{\phi\in R\ ;\
g_\phi(\nu_i)\ge g_\rho(\nu_i)\}$.

In general, the mass of $\rho_I$ on a curve valuation $\nu_\phi$ is the
product of $m(\phi)$ and 
$\div_\phi(I)\=\max\{k\ ;\ \phi^k|\psi\ \text{for all}\ \psi\in I\}$.
%
%
\subsection{Psh functions}\label{S20}
All plurisubharmonic (psh) functions are defined near the origin
in $\C^2$. A psh function $u$ is said to have 
\emph{logarithmic singularities} if there exist $c>0$ and 
holomorphic functions $\phi_i\in R$
such that $u=\frac{c}{2}\log\sum_1^n|\phi_i|^2+O(1)$.
%
%
\subsubsection{Lelong and Kiselman numbers}\label{S10}
For a fixed choice of coordinates $(x,y)$, and weights $a,b>0$ the
\emph{Kiselman number}~\cite{kis1,kis2} of $u$ is defined to be 
\begin{equation*}
  \nu^{x,y}_{a,b}(u)
  =\lim_{r\to0}\frac{ab}{\log r}\sup\{u\ ;\ |x|<r^{1/a},\ |y|<r^{1/b}\}.
\end{equation*}
We have $\nu^{x,y}_{\lambda a,\lambda b}(u)=\lambda\nu^{x,y}_{a,b}(u)$
for any $\lambda>0$. 
When $a=b=1$, the Kiselman number reduces to the 
\emph{Lelong number} $\nu^\mathrm{L}(u)$. The latter does not
depend on the choice of $(x,y)$.
%
%
\begin{figure}[ht]
  \includegraphics[width=0.9\textwidth]{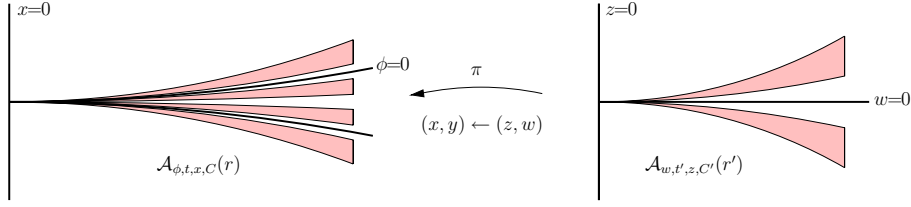}
  \caption{Characteristic regions, see Section~\ref{S2}.
  To the left is a characteristic region $\cA$ around 
  a singular curve $\{\phi=0\}$. If $\pi:X\to(\C^2,0)$ desingularizes
  this curve, the preimage of $\cA$ is a characteristic region
  around a smooth curve.}\label{F2}
\end{figure}
%
%
\subsubsection{Evaluating valuations on psh functions}\label{S2}
(\cite[Section~3]{pshsing}.)
One can evaluate a quasimonomial valuation on a psh function
as the normalized pushforward of a Kiselman number under 
a birational morphism. Concretely, this goes as follows.
If $\nu\in\cVqm$ is quasimonomial, 
write $\nu=\nu_{\phi,t}$ for $\phi\in\fm$ irreducible
with $m(\phi)=m(\nu)=:m$.
Pick a coordinate $x$ transverse to $\phi$, 
\ie $x\in\fm$ is irreducible, $m(x)=1$ and 
$\nu_x\wedge\nu_\phi=\nu_\fm$.
Also pick constants $0<C_1<C_2<\infty$ and consider the region
\begin{equation*}
  \cA(r)
  =\cA_{\phi,t,x,C}(r)
  \=\{|x|<r\ ;\ C_1|x|^{mt}<|\phi|<C_2|x|^{mt}\}
\end{equation*}
for small $r$. See Figure~\ref{F2}.
We have $\vol\cA(r)\simeq r^{2A}$,
where $A=A(\nu)$ is the thinness of $\nu$.
The value of $\nu$ on a psh function $u$ is then given by
\begin{equation*}
  \nu(u)=\lim_{r\to0}\frac1{\log r}\sup_{\cA(r)}u;
\end{equation*}
this does not depend on the choices of $\phi$, $x$ or $C$. 
We always have the upper bound $u(q)\le\nu(u)\log\|q\|+O(1)$
in $\cA(r)$.
If $u=\frac{c}{2}\log\sum_1^n|\phi_i|^2+O(1)$ 
has logarithmic singularities, then
$\nu(u)=c\min_i\nu(\phi_i)=c\,\nu(I)$,
where $I$ is the ideal generated by the $\phi_i$'s. 
In this case we also have the lower
bound $u(q)\ge\nu(u)\log\|q\|+O(1)$ for 
suitable choices of $C_1$ and $C_2$.
See~\cite[Proposition~4.3]{pshsing}.

If $\nu\in\cVqm$ is monomial, \ie of multiplicity one, then
there exist local coordinates $(x,y)$ such that
$\nu=\nu_{y,t}$, where $t=\a(\nu)\ge1$.
In this case $\nu(u)$ equals the Kiselman number
$\nu^{x,y}_{t,1}(u)$.
In particular $\nu_\fm(u)$ is the Lelong number of $u$.
%
%
\subsubsection{Geometric interpretation}\label{S11}
The value of divisorial valuations on psh functions 
can be interpreted along the lines of Section~\ref{S6}.
Using the notation of that section, 
$b(\nu)\nu(u)$ is the Lelong number of the pullback 
$\pi^*u$ at a generic point $p\in E$~\cite[Proposition~4.1]{pshsing}.
Moreover, if $p\in E$ is a smooth point on $\pi^{-1}(0)$,
then the Lelong number of the strict transform of $u$
under $\pi$ at $p$ is bounded from above by 
$b(\nu)^{-1}\rho_u(U(\vv))$~\cite[Lemma~7.6]{pshsing}.
Here $\rho_u$ is the tree measure of $u$ 
(see below) and 
$\vv=\vv_p$ is the tangent vector at $\nu$ associated to $p$,
\ie $\vv$ is represented by any curve valuation $\nu_C$,
such that the strict transform of $C$ is smooth
and intersects $E$ transversely at $p$.
%
%
\subsubsection{Demailly approximation}\label{S1}
Demailly devised a method, crucially important to our analysis,
of approximating a general psh function $u$ with a sequence $u_n$
of psh functions with logarithmic singularities. This goes as
follows: see Theorem~4.2 and its proof in~\cite{DK} for details. 
Suppose $u$ is psh on a fixed ball $B$ containing the origin.
For $n>0$ let $\cH_n$ be the Hilbert space of holomorphic functions 
$\psi$ on $B$ such that $\int_B|\psi|^2e^{-2nu}<\infty$. 
Let $(g_{nl})_1^\infty$ be an orthonormal basis for $\cH_n$ and
set $u_n=\frac1{2n}\log\sum_l|g_{nl}|^2$.
If $0\ni B'\Subset B$ is a smaller ball, then there exists
$k=k(u,n,B')<\infty$ such that 
$u_n-\frac{1}{2n}\log\sum_1^k|g_{nl}|^2$ is bounded in $B'$:
in particular $u_n$ has a logarithmic singularity at the origin.

Demailly approximation interacts very well with the evaluation of
quasimonomial valuations on psh functions: we have
$0\le\nu(u)-\nu(u_n)\le A(\nu)/n$, 
where $A$ denotes thinness~\cite[Proposition~3.12]{pshsing}.
%
%
\subsubsection{Tree transforms}\label{psh-transform}
(\cite[Section~6.1]{pshsing}.)
Any psh function $u$ has an associated \emph{tree transform}
$g_u\in\cP$, defined by $g_u(\nu)=\nu(u)$ for $\nu\in\cVqm$, as well
as an associated \emph{tree measure}
$\rho_u=\Delta g_u$ on $\cV$. 
That these are well defined follows from Demailly approximation. 
If $u$ is psh, then
$\rho_u$ cannot put mass on a formal (\ie non-analytic) curve
valuation. When $\nu_C$ is a curve valuation associated to an
analytic curve $C=\{\phi=0\}$, the mass of $\rho_u$ on $\nu_C$ is exactly
$m(C)$ times the mass of the positive closed $(1,1)$ current $dd^cu$ 
on $C$. By Siu's theorem, this is equivalent to saying that
$u=u'+c\log|\phi|$
with $c=m(C)^{-1}\rho_u\{\nu_\phi\}$ and $u'$ a psh function with
$\rho_{u'}\{\nu_C\}=0$.
%
%
%
%
\section{Multiplier ideals of tree potentials}\label{sec-def}
In this section we define multiplier ideals of tree potentials
and examine their main properties.
As we will show in subsequent sections, this notion contains 
all previously known (to the authors) definitions of multiplier
ideals, in a local, two-dimensional setting. We also
introduce some related singularity exponents.
\begin{Def}\label{our-def}
  Let $h: \cVqm \to [0,\infty)$ be a tree potential on $\cV$, \ie $h\in\cP$.
  We define the \emph{multiplier ideal} 
  $\cJ(h)$ of $h$ to be the ideal of elements $\psi\in R$ 
  such that
  \begin{equation}\label{e801}
    \sup_{\nu\in\cVqm}\frac{h(\nu)}{\nu(\psi)+A(\nu)}<1.
  \end{equation}
\end{Def}
\begin{Def}
  Let $h\in\cP$. The \emph{singularity exponent}, 
  or \emph{log-canonical threshold}, of $h$ 
  is the number 
  $c(h):=\sup\{c>0\ ;\ \cJ(c\,h)=R\}$. 
  The \emph{Arnold multiplicity} is the 
  number $\lambda(h)=c(h)^{-1}$.
  Thus we have
  \begin{equation}\label{e501}
    \lambda(h)=\sup_{\nu\in\cVqm}\frac{h(\nu)}{A(\nu)}
    \qand
    c(h)=\inf_{\nu\in\cVqm}\frac{A(\nu)}{h(\nu)}.
  \end{equation}
\end{Def}
\begin{Remark}\label{R101}
  It follows immediately from the definition that
  $\cJ(c(h)h)\subsetneq R$.
\end{Remark}

It will be useful to have a less elegant but more concrete 
criterion for membership in the multiplier ideal.
\begin{Prop}\label{P301}
  Fix a tree potential $h\in\cP$ and $\psi\in R$. 
  Let $\rho_h$ and $\rho_\psi$ be the associated measures
  on $\cV$.
  Then $\psi\in\cJ(h)$ iff the following hold:
  \begin{itemize}
  \item[(a)]
    $\nu(\psi)+A(\nu)>h(\nu)$ for all quasimonomial valuations $\nu$;
  \item[(b)]
    $\rho_\psi\{\nu\}+m(\nu)>\rho_h\{\nu\}$ for all curve valuations $\nu$.
  \end{itemize}
  Moreover, we have
  \begin{itemize}
  \item[(i)]
    if $\psi\not\in\fm$, then it suffices to check
    (a)-(b) for $\nu$ of multiplicity one;
  \item[(ii)]
    if $h=c\,g_I$ for some ideal $I$ and $c>0$, then it
    suffices to check
    (a)-(b) for $\nu$ in the approximating sequence of 
    some valuation in the support of $\rho_I$.
  \end{itemize}
\end{Prop}
The proof is given below.
Note that~\eqref{e801} is not a trivial consequence of~(a).
%
%
\subsection{Properties}
Let us now state some basic properties of multiplier ideals
associated to tree potentials.
The proofs are given below.
However, two of them are relegated to the appendix 
as they do not follow 
directly from our approach, but rely on the corresponding 
statements for ideals already proved in the literature.
\begin{Prop}\label{P708}
  The multiplier ideal $\cJ(h)$ is integrally closed.
\end{Prop}
It is clear from the definition that if $g\le h$, then
$\cJ(h)\subset\cJ(g)$. 
A deeper result is the following semicontinuity statement.
\begin{Prop}\label{P709}
  If $(h_n)_1^\infty$ is a decreasing sequence of tree potentials
  converging in $\cP$ (\ie pointwise) to a tree potential $h$, 
  then $\cJ(h_n)=\cJ(h)$ for $n\gg1$.
\end{Prop}
Notice that since $R$ is
Noetherian, the sequence $(\cJ(h_n))_1^\infty$ is stationary.
This does not, however, immediately imply the proposition.

We can bound the Arnold multiplicity as follows.
\begin{Prop}\label{P710}
  Consider $h\in\cP$ with associated measure 
  $\rho=\Delta h$. Then:
  \begin{itemize}
  \item[(i)]
    $\frac12 h(\nu_\fm)\le \lambda(h)\le h(\nu_\fm)$;
  \item[(ii)]
    $\lambda(h)>\frac12 h(\nu_\fm)$ iff there exists a tangent vector
    $\vv$ at $\nu_\fm$ with $\rho\,U(\vv)>\frac12$;
  \item[(iii)]
    $\lambda(h)=h(\nu_\fm)$ iff $\rho$ is a point mass at a curve
    valuation of multiplicity one.
  \end{itemize}
\end{Prop}
Next we have the following subadditivity property.
\begin{Prop}\label{P302}
  If $h_1, h_2\in\cP$, then $\cJ(h_1+h_2)\le\cJ(h_1)\cJ(h_2)$.
\end{Prop}
Our proof uses a reduction to the case when the $h_i$'s are
proportional to tree transforms of ideals and is given in 
the appendix.
The same holds for the following version of Skoda's Theorem:
\begin{Prop}\label{P303}
  Let $\nu_0$ be a divisorial valuation with 
  associated (simple complete) ideal
  $I_0=I_{\nu_0}$
  and tree transform $g_0=g_{I_0}$.
  Then we have
  \begin{equation*}
    \cJ(h+2g_0)=I_0\cJ(h+g_0)
  \end{equation*}
  for any tree potential $h\in\cP$.
\end{Prop}
On the other hand, there is a version of Skoda's Theorem
that we can prove directly:
\begin{Prop}~\label{P304}
  Consider $\phi\in\fm$ with associated tree potential $g_\phi$.
  Then
  \begin{equation*}
    \cJ(h+g_\phi)=\phi\cJ(h)
  \end{equation*}
  for any tree potential $h\in\cP$.
\end{Prop}
%
%
\subsection{Proofs}
We now turn to the proofs. The fact that multiplier ideals
are integrally closed is easy to establish.
\begin{proof}[Proof of Proposition~\ref{P708}]
  Let $(\psi_i)_1^k$ be generators of $I=\cJ(h)$. If $\psi\in\bI$,
  then $\nu(\psi)\ge\min_i(\psi_i)$ for every $\nu\in\cVqm$.
  This implies that 
  \begin{equation*}
    \sup_\nu\frac{h(\nu)}{A(\nu)+\nu(\psi)}
    \le\max_i\sup_\nu\frac{h(\nu)}{\nu(\psi_i)+A(\nu)}
    <1
  \end{equation*}
  so that $\psi\in I$.
\end{proof}  

We now turn to the characterization of $\cJ(h)$ given 
in Proposition~\ref{P301}.
The main point is that the supremum in the definition of $\cJ(h)$ 
is actually attained. In fact, the following lemma is the
key to the Openness Conjecture for psh functions.
\begin{Lemma}\label{lem-attain}
  Fix a tree potential $h\in\cP$ and let $\psi\in R$.
  The function
  \begin{equation}\label{def-chi}
    \chi(\nu)=\chi_{h,\psi}(\nu)=\frac{h(\nu)}{\nu(\psi)+A(\nu)},
  \end{equation}
  defined on $\cVqm$, extends to a function on $\cV$ continuous on
  segments. Its supremum is attained at a valuation $\nu_\star$.  We
  may choose $\nu_\star$ to be a quasimonomial or a curve valuation,
  and to belong to the support of $h$. Further, we may take
  $\nu_\star$ such that:
  \begin{itemize}
  \item[(i)]
    if $\psi\notin\fm$, then $\nu_\star$ is of multiplicity one;
  \item[(ii)]
    if $h=g_I$ for some ideal $I$ in $R$, then 
    $\nu_\star$ is an element in 
    the approximating sequence of some valuation in $\supp\rho_I$;
    in particular, $\nu_\star$ is 
    either divisorial or a curve valuation at an element of $R$;
  \item[(iii)]
    if $h=g_u$ for some psh function $u$, then 
    $\nu_\star$ is either quasimonomial 
    or an analytic curve valuation.
  \end{itemize}
\end{Lemma}
The main observation in the proof of Lemma~\ref{lem-attain} is
that it suffices to consider $h$ on a finite subtree of $\cV$.
In the same spirit we have:
\begin{Lemma}\label{L711}
  For $h\in\cP$ and $0<\e<1$, the set
  \begin{equation*}
    \cT
    =\cT_{h,\e}
    :=\{\nu\in\cV\ ;\ 
    \rho_h\{\mu\ge\nu\}\ge(1-\e)m(\nu)\}.
  \end{equation*}
  is a finite subtree of $\cV$ and $\cJ(h)=\cJ(h_\cT)$,
  where $h_\cT$ is the smallest tree potential coinciding with 
  $h$ on $\cT$.
\end{Lemma}
\begin{proof}[Proof of Lemma~\ref{lem-attain}]
  We may assume that $\rho:=\Delta h$ has mass 1.
  Set
  \begin{equation*}
    \cT=\cT_{h,\psi}
      :=\{\nu\ ;\ \rho\{\mu\ge\nu\}\ge \chi(\nu_\fm)\,m(\nu)\}.
  \end{equation*}
  If $\nu\le\nu'$, then $m(\nu)\le m(\nu')$ and
  $\rho\{\mu\ge\nu\}\ge\rho\{\mu\ge\nu'\}$.  Hence $\cT$ is a tree.
  The number of ends of $\cT$ is finite, bounded by 
  $1/\chi(\nu_\fm)$. Thus $\cT$ is a finite tree. Its ends
  have finite multiplicity, hence consists of 
  quasimonomial or curve valuations.
  
  We claim that $\sup_\cV\chi=\sup_\cT\chi$. To see this, pick any
  $\nu_1\not\in\cT$.  Let $\nu_0=\max\{\nu\in\cT\ ;\ \nu\le\nu_1\}$
  and let $\vv$ be the tangent vector at $\nu_0$ represented by
  $\nu_1$. Write $m(\vv)=m$, $\a(\nu_0)=\a_0$ and, for $\nu\in
  [\nu_0,\nu_1]$, $\a=\a(\nu)$.  The right derivative of $h$ at
  $\nu_0$ with respect to the skewness equals 
  $\rho\,U(\vv)$ 
  (see Section~\ref{sec-measure}).  
  As $h$ is concave on $[\nu_0,\nu_1]$ we
  have $h(\nu)\le h(\nu_0)+\rho\,U(\vv)(\a-\a_0)$ for
  $\nu\in[\nu_0,\nu_1]$.  On the other hand, $A(\nu)\ge
  A(\nu_0)+m(\a-\a_0)$ and $\nu(\psi)\ge\nu_0(\psi)$, hence
  \begin{equation*} 
    \chi(\nu)
    \le\frac{h(\nu_0)+\rho\,U(\vv)(\a-\a_0)}
    {\nu_0(\psi)+A(\nu_0)+m(\a-\a_0)} 
    :=M(\a),\quad 
    \text{for all $\nu\in[\nu_0,\nu_1]$}.  
  \end{equation*} 
  Here $\a\mapsto M(\a)$ is M\"obius in $\a$, hence 
  \begin{equation*} 
    \sup_{[\nu_0,\nu_1]}\chi
    \le\sup_{[\a_0,\infty]}M 
    =\max\{\chi(\nu_0),m^{-1}\rho\,U(\vv)\}
    \le\max\{\chi(\nu_0),\chi(\nu_\fm)\}.
  \end{equation*}
  To prove the
  last inequality, we note that $U(\vv)\cap\cT=\emptyset$. 
  If $\rho\,U(\vv)>m\chi(\nu_\fm)$, we could find $\nu'\in U(\vv)$
  close enough to $\nu_0$, of multiplicity $m$ and such 
  that $\rho\{\mu\ge\nu'\}\ge m\chi(\nu_\fm)$. 
  This would imply $\nu'\in\cT$, a contradiction.

  Thus $\sup_\cV\chi=\sup_\cT\chi$.
  
  As each of the functions $h(\nu)$, $\nu(\psi)$ and $A(\nu)$ is
  continuous on segments in $\cVqm$, so is $\chi$. We extend $\chi$ to
  $\cV$ by taking limits along segments.  To see that these limits
  exist, pick $\nu_0\notin\cVqm$ and let $\nu\to\nu_0$ along the
  segment $[\nu_\fm,\nu_0[$. If $\nu_0$ is a curve valuation, then
  $h(\nu)/\a(\nu)\to\rho\{\nu_0\}$,
  $\nu(\psi)/\a(\nu)\to\rho_\psi\{\nu_0\}$ and
  $A(\nu)/\a(\nu)\to m(\nu_0)$. This gives
  $\chi(\nu)\to\rho\{\nu_0\}(\rho_\psi\{\nu_0\}+m(\nu_0))^{-1}$.
  
  When $\nu_0$ is infinitely singular, we claim that $\chi$ is
  decreasing when $\nu\in [\nu_\fm,\nu_0]$ is sufficiently close to
  $\nu_0$. This implies that $\chi(\nu)$ converges when
  $\nu\to\nu_0$. To prove the claim, note that 
  $\nu\mapsto\nu(\psi)$ is constant equal to $\nu_0(\psi)$ near
  $\nu_0$ and that the left derivative of $h$ with respect to 
  skewness at $\nu$ equals $\rho\{\mu\ge\nu\}$. 
  The left derivative of $\chi$ at $\nu$ is hence equal to
  \begin{equation}\label{e001}
    \frac{d\chi}{d\a}
    =\frac{(A(\nu)+\nu_0(\psi))\rho\{\mu\ge\nu\}-m(\nu)h(\nu)}
    {(A(\nu)+\nu_0(\psi))^2}.
  \end{equation}
  By Section~\ref{thinness}, 
  $A(\nu)-m(\nu)\a(\nu)\to-\infty$ as $\nu\to\nu_0$
  and $h(\nu)\ge \rho\{\mu\ge\nu\}\a(\nu)$. 
  This easily implies $d\chi/d\a<0$
  and completes the proof of the claim.
  
  Thus the function $\chi$ extends to $\cV$. 
  Its restriction to the finite tree $\cT$ is continuous, hence 
  $\sup_\cV\chi$ is attained at some $\nu_\star\in\cT$.
  
  Now suppose $\psi\notin\fm$, which amounts to $\nu(\psi)\equiv0$.
  Consider a segment $]\nu_0,\nu_1[$ on which the multiplicity is
  constant equal to $m_0\ge 2$. Then~\eqref{e001} and the two
  inequalities $A(\nu)<m(\nu)\a(\nu)=m_0\a(\nu)$ and
  $h(\nu)\ge\rho\{\mu\ge\nu\}\a(\nu)$, imply $\frac{d\chi}{d\a}<0$ on
  $]\nu_0,\nu_1[$.  This proves that $\chi$ is decreasing off the
  subtree $\{\nu\ ;\ m(\nu)=1\}$. Hence $m(\nu_\star)=1$, proving~(i).
  
  We saw above that if $\nu$ is a curve valuation, then 
  $\chi(\nu)=0$ unless $\rho$ puts mass on $\nu$. This proves~(iii)
  as a measure represented by a psh function cannot charge
  formal curve valuations (see Section~\ref{psh-transform}).

  Finally, if $h=g_I$ for an ideal $I$, then $\rho$ is atomic,
  supported on finitely many valuations that are either divisorial or
  curve valuations, see Section~\ref{S5}. On any open segment in $\cV$
  not intersecting the support of $\rho$, and on which the
  multiplicity is constant, $h$ is M\"obius, hence attains its
  maximum at its boundary points. This implies~(ii).
\end{proof}
\begin{proof}[Proof of Lemma~\ref{L711}]
  Suppose first that $\cT$ is empty. This happens exactly when
  $h(\nu_\fm)<1-\e$. We have $A(\nu) \ge 1 + \a(\nu)$, and $h(\nu)\le
  h(\nu_\fm) \a(\nu)$ for all $\nu$, whence $\sup_\cV\chi_{h,1}\le
  h(\nu_\fm)<1$ so that $\cJ(h)=\cJ(h_\cT)=R$.

  Suppose now that $\cT$ is nonempty \ie $h(\nu_\fm)\ge1-\e$.  We
  follow the proof of Lemma~\ref{lem-attain}. 
  Clearly $\cT$ is a subtree of $\cV$ 
  with at most $h(\nu_\fm)/(1-\e)$ ends. 
  As $h\ge h_\cT$ we have $\cJ(h) \subset \cJ(h_\cT)$. To prove the
  reverse inclusion, we proceed as follows.  Suppose
  $\sup_\cT\chi=t<1$ where $\chi=\chi_{h,\psi}$.  We must show that
  $\sup_\cV\chi<1$.  Consider $\nu_0\in\cT$ and a segment
  $]\nu_0,\nu_1]$, disjoint from $\cT$. Let $m=m(\vv)$, where $\vv$ is
  the tangent vector at $\nu_0$ represented by $\nu_1$.  For
  $\nu\in\,]\nu_0,\nu_1]$ we get
  \begin{equation*}
    \chi(\nu)
    \le\max\{\chi(\nu_0),\rho\{\mu\ge\nu\}/m\}
    \le\max\{t,1-\e\}.
  \end{equation*}
  Hence $\sup_\cV\chi\le\max\{t,1-\e\}<1$,
  which completes the proof.
\end{proof}
\begin{proof}[Proof of Proposition~\ref{P301}]
  This is an immediate consequence of Lemma~\ref{lem-attain}
  and the formula $\chi(\nu)=\rho_h\{\nu\}/(\rho_\psi\{\nu\}+m(\nu))$
  for a curve valuation $\nu$.
\end{proof}

Next we address semicontinuity. The main ingredient 
in the proof is the following lemma. 
It can be viewed as a general statement about
parameterized trees, but we shall formulate it only for the valuative
tree $\cV$ parameterized by skewness.
\begin{Lemma}\label{L708}
  Under the assumptions of Proposition~\ref{P709}, the functions
  $\frac{h_n(\nu)}{\a(\nu)}$ decrease uniformly on $\cVqm$ 
  to the function $h(\nu)/\/a(\nu)$. 
\end{Lemma}
\begin{proof}[Proof of Proposition~\ref{P709}]
  Clearly $\cJ(h_n)$ forms an increasing sequence of ideals in $R$
  and $\cJ(h_n)\subset\cJ(h)$ for all $n$. 
  Since $R$ is Noetherian, it suffices to show that 
  $\bigcup\cJ(h_n)\supset\cJ(h)$. 
  Consider $\phi\in\cJ(h)$ and write
  $\chi(\nu)=h(\nu)/(\nu(\phi)+A(\nu))$ and
  $\chi_n(\nu)=h_n(\nu)/(\nu(\phi)+A(\nu))$.
  Then $\sup\chi<1$. Since $A\ge\a$, Lemma~\ref{L708} gives
  \begin{equation*}
    0\le\chi_n(\nu)-\chi(\nu)
    \le\frac{h_n(\nu)}{\a(\nu)}-\frac{h(\nu)}{\a(\nu)}
    \to0
  \end{equation*}
  uniformly on $\cVqm$, as $n\to\infty$. 
  This implies $\sup\chi_n<1$, \ie $\phi\in\cJ(h_n)$
  for large $n$, completing the proof.
\end{proof}
\begin{proof}[Proof of Lemma~\ref{L708}]
  As mentioned above, a version of the lemma holds for arbitrary 
  parameterized trees. For a totally ordered tree it translates
  into the following elementary statement, the proof of which is
  left to the reader.
  \begin{Lemma}\label{L709}
    Let $g_n$, $g$ be positive concave functions on $[1,\infty)$
    such that $g_n$ decreases to $g$ pointwise. Then 
    $g_n(x)/x$ decreases to $g(x)/x$ uniformly.
  \end{Lemma}
  Of course, a tree such as $\cVqm$ has a lot of branching.
  The following lemma will be used to control the behavior of
  tree potentials at branch points.
  \begin{Lemma}\label{L710}
    Consider tree potentials $g,h$ and quasimonomial valuations 
    $\mu$, $\nu$ with $\mu<\nu$.
    Assume that 
    \begin{equation*}
      \frac{g(\nu)}{\a(\nu)}\ge\frac{h(\nu)}{\a(\nu)}+\e
      \qand
      \frac{g(\mu)}{\a(\mu)}\le\frac{h(\mu)}{\a(\mu)}+\e
    \end{equation*}
    for some $\e>0$. Then $\rho\,U(\vv)\ge\e$, where
    $\rho=\rho_g$ is the tree measure of $g$ and
    $\vv$ is the tangent vector at $\mu$ represented by $\nu$.
  \end{Lemma}
  \begin{proof}
    On the one hand, since $h$ is increasing, the assumptions give
    \begin{equation*}
      g(\nu)
      \ge h(\nu)+\e\a(\nu)
      \ge h(\mu)+\e\a(\nu)
      =h(\mu)+\e\a(\mu)+\e(\a(\nu)-\a(\mu)).
    \end{equation*}
    On the other hand, as $g$ is a tree potential, we have
    \begin{equation*}
      g(\nu)
      \le g(\mu)+\rho\,U(\vv)(\a(\nu)-\a(\mu))
      \le h(\mu)+\e\a(\mu)+\rho\,U(\vv)(\a(\nu)-\a(\mu))
    \end{equation*}
    Thus $\rho\,U(\vv)\ge\e$.
  \end{proof}
  We continue the proof of Lemma~\ref{L708} and argue by 
  contradiction. If $h_n/\a$ does not converge uniformly
  to $h/\a$, then, after passing to a subsequence if necessary,
  there exist $\nu_n\in\cVqm$ and $\e>0$ such that
  $h_n(\nu_n)/\a(\nu_n)>h(\nu_n)/\a(\nu_n)+\e$ for all $n$. 
  Using the compactness of $\cV$ we may assume that 
  $\nu_n$ converges weakly to some valuation $\nu_\star\in\cV$.
  Since $h_n$ converges to $h$ pointwise on $\cVqm$ we may
  assume that the $\nu_n$ are all distinct and different from
  $\nu_\star$.

  By Lemma~\ref{L709} we may assume 
  $h_n/\a\le h/\a+\e$ on $[\nu_\fm,\nu_\star[$. 
  Hence $\nu_n\notin[\nu_\fm,\nu_\star]$ 
  so that $\mu_n:=\nu_\star\wedge\nu_n$ is quasimonomial 
  and $\mu_n<\nu_n$. We apply Lemma~\ref{L710} 
  to $g=h_n$, $h=h$, $\mu=\mu_k$ and $\nu=\nu_k$, for 
  $n\ge k$. Recall that $h_n\ge h_k$. We conclude that
  $\rho_n\,U_k\ge\e$, where $\rho_n=\rho_{h_n}$ and
  $U_k=U(\vv_k)$, where $\vv_k$ is the tangent vector
  at $\mu_k$ represented by $\nu_k$.

  If $k\ne l$, then $U_k=U_l$ or $U_k\cap U_l=\emptyset$.
  As $\mass\rho_n=h_n(\nu_\fm)\to h(\nu_\fm)$ is uniformly
  bounded, this implies that, after passing to a subsequence,
  all the $U_k$ are equal. In particular all the $\mu_k$ are equal.
  As $\mu_k\to\nu_\star$, this gives 
  $\mu_k=\nu_\star$ for all $k$. In particular, $\nu_\star$ is
  quasimonomial and $\nu_n>\nu_\star$.
  
  For $k>1$ we have $\nu_k\wedge\nu_1\in\,]\nu_\star,\nu_1]$ as
  $U_k=U_1$. Again by Lemma~\ref{L709}, $\nu_k\wedge \nu_1 <\nu_k$ for all
  $k\ge2$. As $\nu_1\wedge\nu_k\to\nu_\star$ as $k\to\infty$ we may
  assume that the valuations $\nu_1\wedge\nu_k$ are all distinct. We
  again apply Lemma~\ref{L710}, now with $\mu=\nu_1\wedge\nu_k$,
  $\nu=\nu_k$ and still using $h=h$ and $g=h_n$ for $n\ge k$. We
  conclude that $\rho_n\,U'_k\ge\e$, where $U'_k=U(\vv'_k)$ and
  $\vv'_k$ is the tangent vector at $\nu_1\wedge\nu_k$ represented by
  $\nu_k$. But the regions $U'_k$ are mutually disjoint, so since the
  $\rho_n$ has uniformly bounded mass, this gives a contradiction.
\end{proof}

Next we prove the bounds for the singularity exponents.
\begin{proof}[Proof of Proposition~\ref{P710}]
  Recall that $\rho$ has mass $h(\nu_\fm)$.

  Write $\chi(\nu)=h(\nu)/A(\nu)$ so that $\lambda(h)=\sup\chi$.
  The lower bound $\lambda(h)\ge\chi(\nu_\fm)=\frac12 h(\nu_\fm)$
  is immediate.
  Recall that $A(\nu)\ge 1+\a(\nu)$.
  If $\nu\in\cVqm$, $\nu\ne\nu_\fm$, then 
  $h(\nu)\le h(\nu_\fm)+\rho\,U(\vv)(\a(\nu)-1)$,
  where $\vv$ is the tangent vector at $\nu_\fm$ represented
  by $\nu$. 
  In general, this gives $h(\nu)\le h(\nu_\fm)\a(\nu)$, so
  that $\chi(\nu)<h(\nu_\fm)$.

  If $\rho\,U(\vv)\le\frac12 h(\nu_\fm)$, then 
  $h(\nu)\le\frac12h(\nu_\fm)(1+\a(\nu))$ so that
  $\chi(\nu)\le\frac12 h(\nu_\fm)$. 
  On the other hand, if 
  $\rho\,U(\vv)>\frac12 h(\nu_\fm)$, 
  then we may move $\nu$ closer to $\nu_\fm$, still
  keeping $\nu\in U(\vv)$, 
  so that $m(\nu)=1$, $A(\nu)=1+\a(\nu)$ 
  and $h(\nu)>h(\nu_\fm)+\frac12h(\nu_\fm)(\a(\nu)-1)$. This
  gives $\lambda(h)\ge\chi(\nu)>\frac12 h(\nu_\fm)$. 

  Finally suppose $\lambda(h)=h(\nu_\fm)$.
  By Lemma~\ref{lem-attain}
  we have $\lambda(h)=\chi(\nu)$ for some valuation $\nu\in\cV$ 
  of multiplicity one. We have seen above that $\chi< h(\nu_\fm)$ 
  on $\cVqm$. Hence $\nu$ is a curve valuation. But then
  $\chi(\nu)=\rho\{\nu\}/m(\nu)=\rho\{\nu\}$ so that 
  $\rho\{\nu\}=h(\nu_\fm)$ and
  $\rho$ is a point mass at $\nu$.
\end{proof}

Finally we prove the second version of Skoda's Theorem.
\begin{proof}[Proof of Proposition~\ref{P304}]
  By unique factorization we may assume that 
  $\phi$ is irreducible.
  We first claim that $J(h+g_\phi)\subset\phi R$.
  Indeed, if $\psi\in\cJ(h+g_\phi)$, then by
  Proposition~\ref{P301} we have 
  $\rho_\psi\{\nu_\phi\}+m(\phi)>\rho_h\{\nu_\phi\}+\rho_\phi\{\nu_\phi\}
  \ge0+m(\phi)$. Thus $\rho_\psi\{\nu_\phi\}>0$ so that 
  $\phi$ divides $\psi$.

  It hence suffices to show that if $\psi\in R$, then 
  $\phi\psi\in\cJ(h+g_\phi)$ iff $\psi\in\cJ(h)$.
  But this follows easily from Proposition~\ref{P301}. Indeed,
  if $\nu$ is quasimonomial, then 
  $\nu(\psi\phi)-(h+g_\phi)(\nu)=\nu(\psi)-h(\nu)$
  and if $\nu$ is a curve valuation, then 
  $\rho_{\psi\phi}\{\nu\}-(\rho_h+\rho_\phi)\{\nu\}
  =\rho_\psi\{\nu\}-\rho_h\{\nu\}$.
\end{proof}
%
%
%
%
\section{Multiplier ideals of ideals}\label{sec-ideal}
We now show that if $h$ is the tree transform of 
a formal power $I^c$, where $I$ is an ideal and $c>0$,
then the multiplier ideal of 
$h$ agrees with multiplier ideal of $I^c$ 
with the standard definitions given in the literature. 
First we consider the approach based on 
resolution of singularities 
as described in~\cite[Part~III]{laz-book}.
Then we explain how to recover Lipman's definition given in~\cite{Lip94}.
Finally we show how to understand asymptotic multiplier ideals
of graded systems of ideals in terms of tree potentials.
%
%
\subsection{Multiplier ideals through log-resolutions}\label{S301}
The exposition follows~\cite[Part~III]{laz-book}.
See also~\cite{LipWat03} and~\cite{TakWat03}.

A \emph{log-resolution} of an ideal $I\subset R$ is a
birational morphism $\pi:X\to\Spec R$ such that $X$ is non-singular,
$I\cdot\cO_X$ is locally principal: $I\cdot\cO_X=\cO_X(-F)$
for an effective divisor $F$, and
$F+E$ has simple normal crossing singularities,
where $E$ is the exceptional divisor of $\pi$.

The existence of such a resolution in our setting is 
well-known. In fact, $\pi$ can be chosen as a composition
of blowups at (closed) points. 
In the analytic case $R=\cO_0$ the reader may think of 
$\Spec R$ as a neighborhood of the origin in $\C^2$.

Another ingredient in the definition is 
\emph{the relative canonical divisor} of $\pi$. 
It is the unique divisor $K_{X/R}$, with support in
$E$ and whose class represents $J\pi$,
the Jacobian determinant of $\pi$
\begin{equation*}
  K_{X/R}=K_X-\pi^*K_{\Spec R}.
\end{equation*}
Here $K_X$ and $K_{\Spec R}$ denote 
the sheaf of $2$-differentials over $X$ and $\Spec R$, 
respectively.  

We also denote by $\lfloor c \rfloor$ the round-down of a positive
real number \ie the greatest integer less than or equal to $c$.
\begin{Def}\label{classical-def}
  Let $I\subset R$ be an ideal, and $c$ be a positive real number.  Fix
  a log-resolution $\pi:X\to\Spec R$ of $I$, such that 
  $I\cdot\cO_X=\cO_X(-F)$ with $F=\sum r_iE_i$, and write 
  $K_{X/R}=\sum\b_iE_i$.
  We define \emph{multiplier ideal} of $I^c$ to be
  \begin{equation}\label{eq-def-mul}
    \cJ(I^c) 
    \=\{\psi\in R\ ;\ \div_{E_i}(\pi^*\psi)
    \ge\lfloor cr_i\rfloor-\b_i 
    \text{ for all $i$}\}.\tag{$\star$}
  \end{equation}
  Here $\div_{E_i}$ denotes the generic order of vanishing 
  along $E_i$ of a function on $X$.
\end{Def}
This definition does not depend on the choice of log-resolution.
In the analytic case $R=\cO_0$, the definition says that 
a holomorphic germ $\psi$ belongs to the multiplier ideal
iff the pullback of $\psi$ vanishes to sufficiently high orders 
along the irreducible 
components of $\pi^{-1}(V)$, where $V=\bigcap_{\phi\in I}\phi^{-1}(0)$.

\begin{Prop}\label{P401}
  For any ideal $I\subset R$ 
  and any positive real number
  $c>0$, the multiplier ideal $\cJ(I^c)$ as defined
  in Definition~\ref{classical-def} coincides with the multiplier ideal 
  of the tree transform $c\,g_I$ as defined in
  Definition~\ref{our-def}.
\end{Prop}
\begin{Remark}
  Similarly, we may associate a multiplier ideal to an expression
  of the form $I_1^{c_1}\cdots I_n^{c_n}$ for ideals $I_j$ and
  real numbers $c_j$. The proof of Proposition~\ref{P401} is easily
  adapted to show that
  $J(I_1^{c_1}\cdots I_n^{c_n})=J(\sum_jc_j\,g_{I_j})$.
\end{Remark}
\begin{proof}
  We first need to translate the condition of vanishing of $\pi^*\psi$
  along $E_i\subset X$ in terms of valuations. We shall denote this
  condition by $(\star)_i$. Note that $\div_{E_i}$ is a
  valuation with values in $\Z$. 
  We have $r_i=\div_{E_i}(\pi^*I)$, and $\b_i=\div_{E_i}(J\pi)$ 
  where $J\pi$ is the Jacobian determinant of $\pi$.
  
  First suppose $E_i$ does not belong to the exceptional divisor of
  $\pi$. Its image by $\pi$ is then an irreducible curve.
  We let $\nu_i\in\cV$ be the associated curve valuation
  (see Section~\ref{curve-val}) and $m_i$ its multiplicity.
  Denote by $\rho_I$ and $\rho_\psi$ the tree measures 
  of $I$ and $\psi$, respectively.
  It follows from Section~\ref{S5} that $\div_{E_i}(\pi^*I)$ is equal to
  $m_i^{-1}\,\rho_I\{\nu_i\}$. 
  In the same way,
  $\div_{E_i}(\pi^*\psi)=m_i^{-1}\,\rho_\psi\{\nu_i\}$.  
  On the other hand, it is clear that $\b_i=\div_{E_i}(J\pi)=0$.
  The condition~\eqref{eq-def-mul}$_i$ is hence equivalent to
  $m_i^{-1}\rho_\psi\{\nu_i\}\ge\lfloor cm_i^{-1}\rho_I\{\nu_i\}\rfloor$. 
  As the left hand side is an integer, this amounts to
  $m_i^{-1}\rho_\psi\{\nu_i\}>c m_i^{-1}\rho_I\{\nu_i\}-1$. 
  We thus conclude that if $E_i\not\subset\pi^{-1}(0)$, then 
  \begin{equation}\label{E301}
    (\star)_i 
    \Leftrightarrow 
    \rho_\psi\{\nu_i\}+m(\nu_i)>c\rho_I\{\nu_i\}
  \end{equation}

  Now suppose $E_i$ is an irreducible component of $\pi^{-1}(0)$.
  Then $\pi_*\div_{E_i}(\psi)\=\div_{E_i}(\pi^*\psi)$ defines 
  a divisorial valuation with values in $\Z$. 
  It is proportional to a unique normalized valuation $\nu_i\in\cV$: 
  $\pi_*\div_{E_i}=b_i\nu_i$
  where $b_i=b(\nu_i)$ is the generic multiplicity of $\nu_i$ 
  (see Section~\ref{S6}). 
  The thinness $A(\nu_i)$ is given by $a_i/b_i$, where 
  $a_i=a(\nu_i)=\div_{E_i}(J\pi)+1$ (again by Section~\ref{S6}). 
  We infer that if $E_i\subset\pi^{-1}(0)$, then 
  \begin{align*}
    (\star)_i
    \Leftrightarrow\ 
    &\div_{E_i}(\pi^*\psi)
    \ge\lfloor c\,\div_{E_i}(\pi^*I)\rfloor-\div{E_i}(J\pi)\\
    \Leftrightarrow\ 
    &\div_{E_i}(\pi^*\psi)>c\,\div_{E_i}(\pi^*I)-\div{E_i}(J\pi)-1\\
    \Leftrightarrow\ 
    &b_i\nu_i(\psi)>c\,b_i\nu_i(I)-a_i
  \end{align*}
  \ie
  \begin{equation}\label{E401}
    (\star)_i
    \Leftrightarrow
    \nu_i(\psi)+A(\nu_i)>c\,g_I(\nu_i).
  \end{equation}

  Now consider $\psi\in\cJ(c\,g_I)$. It follows immediately
  from Proposition~\ref{P401} that 
  condition $(\star)_i$ in~\eqref{E301},~\eqref{E401} holds for
  all $i$. Thus $\psi\in\cJ(I^c)$.

  Conversely suppose $\psi\not\in\cJ(c\,g_I)$. Again we use
  Proposition~\ref{P301}. This tells us that either~(a) or~(b) fails 
  in that proposition. 
  If~(b) fails, then there exists a curve valuation
  $\nu$ such that 
  $\rho_\psi\{\nu\}+m(\nu)\le c\rho_I\{\nu\}$. 
  This inequality implies that the element in $R$ corresponding to
  $\nu$ is a factor of $I$. Hence $\nu$ is one of the curve valuations
  $\nu_i$ above. Thus $(\star)_i$ in~\eqref{E301} fails and 
  $\psi\not\in\cJ(I^c)$.
  On the other hand, if~(a) fails, then there exists
  a divisorial valuation $\nu$ such that 
  $\nu(\psi)+A(\nu)\le c\,\nu(I)$. 
  We may take a log-resolution $\pi$ of $I$ such that 
  $\nu$ corresponds to an exceptional divisor of $\pi$.
  This can be achieved by further blowups of a given
  log-resolution of $I$. Thus $\nu=\nu_i$ for some $i$,
  so that~$(\star)_i$ in~\eqref{E401} fails and $\psi\not\in\cJ(I^c)$.
\end{proof}
%
%
\subsection{Lipman's approach}\label{sec-lipmanapproach}
The definition of 
Lipman\footnote{Lipman used the term \emph{adjoint ideals} 
  and defined them for ideals rather than formal powers of ideals.}
in~\cite{Lip94}
has the advantage of not depending on a 
log-resolution and hence makes sense in great generality.

We describe Lipman's construction in our setting, namely for
an equicharacteristic zero, two-dimensional, regular local ring $R$
with algebraically closed residue field $k$.
Write $K$ for the fraction field of $R$.

We first need some definitions. For a (not necessarily centered)
valuation $\mu$ on $R$, 
we let $R_\mu=\{\phi\in K\ ;\ \mu(\phi)\ge0\}$ be the valuation ring, 
$\fm_\mu=\{\phi\in K\ ;\ \mu(\phi)>0\}$ the (unique) maximal ideal 
in $R_\mu$ and $k_\mu:=R_\mu/\fm_\mu$ the residue field.
Note that  $k_\mu$ is a field extension of $k$.
A \emph{prime divisor} is by definition a valuation on $R$
whose residue field has transcendence degree $l-1$ over $k$, where
$l$ is the height of $\fm_\mu\cap R$ \ie the dimension of the center
of $\mu$ in $\Spec R$.  Let us describe explicitly in our situation
what a prime divisor means. Two different cases may appear.
\begin{itemize}
\item[(i)]
  $l=1$. Then the center of $\mu$ is a principal ideal
  generated by an irreducible element $\psi\in R$. 
  In this case for any $\phi\in R$, 
  we have $\mu(\phi)=\div_C(\phi)$ for some irreducible
  $C\subset\Spec R$.
\item[(ii)] 
  $l=2$. Then $\mu$ is centered at the maximal ideal $\fm$ of $R$.
  The transcendence degree of $k_\mu$ being $1$ implies
  that $\mu$ is a divisorial valuation 
  (see~\cite[Proposition 1.12]{valtree}).
\end{itemize}
Lipman's construction uses the \emph{Jacobian ideal} $J_{R_\mu/R}$,
whose definition is purely algebraic. In our setting it is given as
follows, following the two cases above:
\begin{itemize}
\item[(i)]
  When $l=1$, the ideal $J_{R_\mu/R}$ is trivial. 
\item[(ii)]  
  When $l=2$, $\mu$ is a divisorial valuation. 
  Fix a composition of blowups $\pi:X\to\Spec R$ such that
  the center of $\mu$ in $X$ is a divisor $E$. 
  For a closed point $p\in E$, the ring $\mathcal{O}_{X,p}$ 
  can be naturally viewed as a subring
  of $R_\mu$ by the isomorphism $\pi_*$ between function fields of $X$
  and $\Spec R$.  Then the ideal $J_{R_\mu/R}\subset R_\mu$ is
  generated by the Jacobian determinant of $\pi$.
\end{itemize}
We can now set
\begin{Def}\label{lipman-def}
  Let $I\subset R$ be an ideal and let $c>0$. 
  The multiplier ideal of $I^c$ is the ideal
 \begin{equation*}
   \cJ(I^c)\=\bigcap_\mu\{\psi\in R\ ;\ 
   \mu(\psi)\ge\lfloor c\mu(I)\rfloor-\mu(J_{R_\mu/R})\}
 \end{equation*}
 where the intersection is taken over all prime divisors $\mu$ 
 of $R$.\footnote{Lipman's original definition was for 
   $c=1$ in which case the ``$\lfloor\cdot\rfloor$'' can be omitted.}
\end{Def}
If $I$ is an ideal and $c>0$, then we may consider both the multiplier
ideal of $I^c$ above and the multiplier ideal of the tree potential
$c\,g_I$. Both are defined in terms of valuations. Hence the following
result is perhaps not so surprising.
\begin{Prop}
  For any ideal $I\subset R$ and any $c>0$, 
  the multiplier ideal $\cJ(I^c)$ in the 
  sense of Definition~\ref{lipman-def} coincides with the multiplier 
  ideal of the tree potential $c\,g_I$ as defined in 
  Definition~\ref{our-def}.
\end{Prop}
\begin{proof}
  Pick $\psi\in R$.
  First consider a curve valuation $\nu=\nu_C$ in $\cV$.
  We denote the prime divisor associated to $\nu$ by $\div_C$. 
  Note that $\div_C$ is not centered at the maximal ideal $\fm\subset R$. 
  We then have the  sequence of equivalences
  \begin{align*}
    \div_C(\psi)\ge\lfloor c\div_C(I)\rfloor-\div_C (J_{R_{\div_C}/R})
    \Leftrightarrow\ &
    \div_C(\psi) > c\div_C(I)-1
    \\ \Leftrightarrow\ & 
    \rho_\psi\{\nu_C\}> c\rho_I\{\nu_C\}-m(C).
  \end{align*}

  Now consider a divisorial valuation $\nu\in\cV$.
  Fix a composition of blowups 
  $\pi$ such that $b\,\nu=\pi_*\div_E$ for some exceptional
  component $E$, $b$ being the generic multiplicity of $\nu$ 
  (see Section~\ref{S6}).
  Let $a-1$ be the order of
  vanishing of the Jacobian determinant of $\pi$ along $E$. 
  By the discussion above
  \begin{align*}
    \nu(\psi)\ge\lfloor c\nu(I)\rfloor-\nu(J_{R_\nu/R})
    \Leftrightarrow\ &
    \div_E(\pi^*\psi)\ge\lfloor c\div_E(\pi^*I)\rfloor-\div_E(J\pi)
    \\ \Leftrightarrow\ &
    \div_E(\pi^*\psi)>c\div_E(\pi^*I) -\div_E(J\pi)-1
    \\ \Leftrightarrow\ &
    b\,\nu(\psi)>bc\,\nu(I)-a
    \\ \Leftrightarrow\ &
    \nu(\psi)>c\,g_I(\nu)-A(\nu).
  \end{align*}
  This concludes the proof in view of Proposition~\ref{P301}.
\end{proof}
%
%
\subsection{Multiplier ideals of graded systems}\label{sec-grad}
Another type of multiplier ideals have been defined for 
graded systems of ideals. Here we show that these multiplier ideals
can also be understood through tree potentials.

Let $I_\bullet=(I_k)_{k\ge1}$ be a graded system of ideals. 
This means that $I_k$ is an ideal in $R$ for every $k$ and that 
$I_kI_l\subset I_{k+l}$ for all $k,l$. Fix $c>0$ rational.
It is proved in~\cite{laz-book} (see also below) 
that the sequence of ideals $(\cJ(I_k^{c/k}))_{k\ge1}$ has a unique 
maximal element.
\begin{Def}
  The \emph{asymptotic multiplier ideal} $\cJ(I_\bullet^c)$ is the unique
  maximal element of the sequence $(\cJ(I_k^{c/k}))_{k\ge1}$.
\end{Def}
\begin{Def}
  The \emph{singularity exponent} 
  (or \emph{log-canonical threshold}) 
  of $I_\bullet$ is
  $c(I_\bullet):=\sup\{c>0\ ;\ \cJ(I_\bullet^c)=R\}$.
  The \emph{Arnold multiplicity} is 
  $\lambda(I_\bullet):=c(I_\bullet)^{-1}$.
\end{Def}
Let $g_k$ be the tree transform of $I_k$. The condition 
$I_kI_l\subset I_{k+l}$ then implies $g_{k+l}\le g_k+g_l$.
It is then elementary that the sequence $k^{-1}g_k$ converges
in $\cP$ (\ie pointwise) to the tree potential
$g:=\inf k^{-1}g_k$. Moreover, the subsequence 
$2^{-j}g_{2^j}$ is decreasing, so by semicontinuity 
(Proposition~\ref{P709})
we have $\cJ(c2^{-j}g_{2^j})=\cJ(c\,h)$ for $j\gg0$. On the
other hand, $\cJ(ck^{-1}g_k)=\cJ(I_k^{c/k})$ for all $k$
by Proposition~\ref{P401}. This gives
\begin{Prop}\label{P305}
  Given a graded sequence $I_\bullet$ of ideals
  there exists a tree potential $g=g_{I_\bullet}$ such that
  the asymptotic multiplier ideal $\cJ(I_\bullet^c)$ 
  coincides with the multiplier ideal $\cJ(c\,g)$
  for any $c>0$.
  As a consequence, $c(I_\bullet)=c(g)$ 
  and $\lambda(I_\bullet)=\lambda(g)$.
\end{Prop}
\begin{Remark}
  It follows from Theorem~\ref{thm-apx-gen} 
  that $g=g_{I_\bullet}$ is uniquely determined by 
  the properties stated in Proposition~\ref{P305}.
\end{Remark}
\begin{Remark}
  Not every tree potential $g$ is associated to 
  a graded system of ideals $I_\bullet$. An example is given by
  $g=\a(\mu\wedge\cdot)$, where $\mu$ is infinitely singular of
  infinite skewness. Indeed, $g(\mu)=\infty$, whereas 
  $\mu(I)<\infty$ for every ideal $I$, 
  hence $g_{I_\bullet}(\mu)<\infty$ for any graded system $I_\bullet$.
  It would be interesting to characterize 
  all tree potentials associated to graded systems of ideals.
  See~\cite[Section~6.2]{pshsing} for related questions.
\end{Remark}
\begin{Example}\label{E502}
  If $I\subset R$ is a fixed ideal and $k_0\in\N$, then
  $I_k:=I^{k+k_0}$ defines a graded system of ideals.
  In this case $k^{-1}g_{I_k}=(1+k_0/k)g_I\to g_I$ so 
  $\cJ(I_\bullet^c)=\cJ(I^c)$ independently of $k_0$.
\end{Example}
\begin{Example}\label{E501}
  If $\nu\in\cV$, then $I_k:=\{\phi\in R\ ;\ \nu(\phi)\ge k\}$ defines
  a graded system of ideals. In this case the tree potential $h$ in
  Proposition~\ref{P305} satisfies
  \begin{equation}\label{e000} 
    \Delta h=\a(\nu)^{-1}\nu.
  \end{equation}
  Moreover,
  $c(I_\bullet)=\a(\nu)(1+\a(\nu_1)^{-1})$, where $\nu_1$ is the first
  element in the approximating sequence of $\nu$.  In particular, if
  $\a(\nu)=\infty$, then $g=0$, and $c(I_\bullet)=\infty$.

  Let us sketch a proof of~\eqref{e000}. It is easy to see in general
  that $g(\mu)\ge h(\mu)\=\a(\nu)^{-1}\a(\mu\wedge\nu)$. For the
  reverse inequality, suppose first that $\nu$ is a divisorial
  valuation of generic multiplicity $b$. Then the tree transform of
  $I_{b\a(\nu)}$ is precisely $\mu\to b\a(\mu\wedge\nu)$,
  see Section~\ref{S5}. Whence $g=h$ in this case. In general, take an
  increasing sequence of divisorial valuations $\nu_n$ tending to $\nu$,
  and apply the preceding result.
\end{Example}
%
%
%
%
\section{Multiplier ideals of psh functions}\label{sec-psh}
We now turn to multiplier ideals of psh functions. For this to make
sense we work in the analytic case $R=\cO_0$. 
Recall that if $u$ is psh, then the 
\emph{multiplier ideal} $\cJ(u)$ consists of the holomorphic
germs $\psi\in R$ such that $|\psi|^2e^{-2u}$ is 
locally integrable at the origin (see~\cite{nadel,DK}). 
The (complex) \emph{singularity exponent} 
(or \emph{log-canonical threshold}) of $u$
at the origin is $c(u)=\sup\{c>0\ ;\ e^{-2cu}\in\Lloc\}$.
The \emph{Arnold multiplicity} is $\lambda(u)=c(u)^{-1}$.
Our goal is to prove
\begin{Thm}\label{T703}
  If $u$ is a psh function then the multiplier ideal
  of $u$ equals that of its tree transform. In other words,
  if $\psi$ is a holomorphic germ, then
  \begin{equation}\label{e709}
    \psi\in\cJ(u)
    \quad\text{iff}\quad
    \sup_{\nu\in\cVqm}\frac{\nu(u)}{\nu(\psi)+A(\nu)}<1.
  \end{equation}
\end{Thm}
The proof makes essential use of the Demailly approximation
technique in Section~\ref{S1}. In fact we have
\begin{Cor}
  If $u_n$ are the Demailly approximants of $u$, then 
  $\cJ(u_n)=\cJ(u)$ for $n\gg1$.
\end{Cor}
\begin{proof}
  We have $0\le\nu(u)-\nu(u_n)\le A(\nu)/n$ for all 
  $\nu\in\cVqm$, hence
  \begin{equation*}
    0
    \le\frac{\nu(u)}{\nu(\psi)+A(\nu)}-\frac{\nu(u_n)}{\nu(\psi)+A(\nu)}
    \le\frac1n,
  \end{equation*}
  which easily implies the corollary in view of~\eqref{e709}.
\end{proof}
We split the proof of the theorem into two parts: 
non-integrability and integrability. 
Set
\begin{equation*}
  \chi(\nu)
  =\chi_{u,\psi}(\nu)
  =\frac{\nu(u)}{\nu(\psi)+A(\nu)}.
\end{equation*}
By Lemma~\ref{lem-attain}, $\chi$ attains it supremum
on $\cV$. We emphasize that this fact is really the key to the 
openness conjecture (Corollary~\ref{cor1}). 
At any rate, Theorem~\ref{T703} follows immediately from
Proposition~\ref{prop-nonint} and 
Proposition~\ref{prop-int} below. 
\begin{Prop}\label{prop-nonint}
  Let $u$ and $\psi$ be as in Theorem~\ref{T703}.
  Suppose there exists $\nu\in\cV$ such that
  $\chi(\nu)\ge1$. Then $|\psi|^2e^{-2u}$ 
  is not locally integrable at the origin.
  In fact, we then have $\vol\{|\psi|^2e^{-2u}>R\}\gtrsim R^{-1}$
  as $R\to\infty$.
\end{Prop}
\begin{Prop}\label{prop-int}
  Let $u$ and $\psi$ be as in Theorem~\ref{T703}.
  Suppose that $\sup\chi<1$.
  Then $|\psi|^2e^{-2u}$ is locally integrable at the origin.
\end{Prop}
\begin{proof}[Proof of Proposition~\ref{prop-nonint}]
  It suffices to show the estimate
  $\vol\{|\psi|^2e^{-2u}>R\}\gtrsim R^{-1}$ as $R\to\infty$.
  By Lemma~\ref{lem-attain}, we may assume $\nu$ is not
  infinitely singular. 

  First suppose $\nu$ is quasimonomial.  
  Then $\nu(\psi)-\nu(u)\le-A(\nu)<0$. 
  We use the analysis in Section~\ref{S2}.
  Write $\nu=\nu_{\phi,t}$ and pick a coordinate $x$ 
  transverse $\phi$.
  Let $\cA=\cA_{\phi,t,x,C}(r)$ be a characteristic region.
  Then $\vol(\cA)\gtrsim r^{2A(\nu)}$.
  Moreover, $u(p)\le\nu(u)\log\|p\|r+O(1)$ and
  $\log|\psi|\ge \nu(\psi)\log\|p\|+O(1)$, so
  \begin{equation*}
    \log|\psi|-u
    \ge(\nu(\psi)-\nu(u))\log\|p\|+O(1)    
    \ge-A(\nu)\log r+O(1)
  \end{equation*}
  in $\cA$.
  The desired estimate follows by choosing $R\sim r^{-2A(\nu)}$.

  Now suppose $\nu=\nu_\phi$ is a curve valuation. Then $\phi$
  is a holomorphic germ. By Section~\ref{psh-transform}, 
  we have $u=c\log|\phi|+u'$ with $c=\rho_u\{\nu_\phi\}/m(\phi)$
  and $u'$ psh.
  Similarly, $\psi=\phi^a\psi'$, where $a\ge0$ and
  $\phi$ does not divide $\psi'$. 
  Then $1\le\chi(\nu_\phi)=c/(a+1)$, so that $c\ge a+1$.
  Since $u'$ is bounded from above we have
  \begin{equation*}
    |\psi|^2e^{-2u}
    \gtrsim|\phi|^{2(a-c)}|\psi'|^2
    \gtrsim\frac{|\psi'|^2}{|\phi|^{2}}
  \end{equation*}
  near the origin. 
  Pick a composition of blowups $\pi$, and a smooth 
  point $p\in\pi^{-1}(0)$ 
  such that the strict transform of $\phi^{-1}(0)$ is
  smooth and intersects the exceptional divisor transversely at $p$,
  and such that the strict transform of $(\psi')^{-1}(0)$ does not 
  contain $p$. This can be done by successively blowing up
  the intersection point of the strict transform of 
  $\phi^{-1}(0)$ and the exceptional divisor.
  Fix coordinates $(z,w)$ at $p$ such that $\pi^{-1}(0)=\{z=0\}$ and
  $\phi^{-1}(0)=\{zw=0\}$. 
  We may assume $\pi^*\phi=z^kw$, $J\pi=z^l$ and
  $\pi^*\psi'=z^m\tilde\psi$, where $k,l,m\in\N^*$
  and $\tilde\psi(0)\ne0$.
  Fix $r_0>0$ small, and consider $0<r\ll r_0$. 
  Set $\Omega(r):=\pi\{|z|<r_0, |w|<r\}$.
  By the change of variables formula, $\vol\Omega(r)\gtrsim r^2$.
  If $p=\pi(z,w)\in\Omega(r)$, then
  \begin{equation*}
    |\psi(p)|^2e^{-2u(p)}
    \gtrsim\frac{|\pi^*\psi'(z,w)|^2}{|\pi^*\phi(z,w)|^{2}}
    \gtrsim|w|^{-2}
    \ge r^{-2}.
  \end{equation*}
  The desired estimate follows by choosing $R\sim r^{-2}$.
\end{proof}
\begin{proof}[Proof of Proposition~\ref{prop-int}]
  First suppose $u$ has logarithmic singularities, 
  say $u=\frac{c}2\log\sum_{l=1}^k|\phi_l|^2$ for $\phi_i\in R$
  and $c>0$. The tree transform of $u$ then coincides with
  the tree transform of $I^c$, where $I$ is the ideal generated
  by the $\phi_i$'s. Integrability can now be proved as
  in~\cite[Proposition~1.7]{DK}, using the results of
  Section~\ref{sec-ideal}.
  
  This goes as follows. 
  Let $\pi:X\to(\C^2,0)$ be a log-resolution of the ideal $I$.
  Thus the total transform of the curve 
  $V=\{\prod_l\phi_l=0\}$ is a union of smooth components with 
  normal crossing singularities. 
  By the change of variables formula, the function 
  $\Phi^2=|\psi|^2\exp(-2u)$ is locally integrable iff 
  $(\pi^*\Phi\cdot|J\pi|)^2$ is locally integrable at 
  any point $p\in\pi^{-1}(0)$.

  Let $\{E_i\}$ be the set of irreducible components of 
  $\pi^{-1}(V)$. For any $i$, we let $\beta_i$, $r_i$ and $\gamma_i$
  be the order of vanishing along $E_i$ of 
  $J\pi$, $\pi^*I$ and $\pi^*\psi$, respectively.
  As in Section~\ref{sec-ideal} the condition
  $\sup\chi<1$ implies that
  $\delta_i:=\gamma_i+\beta_i-\lfloor cr_i\rfloor\ge0$,
  hence $\delta_i>-1$ for all $i$.

  Now pick $p\in\pi^{-1}(0)$ and local
  coordinates $(z,w)$ at $p$ such that 
  $\pi^{-1}V\subset\{zw=0\}$. The calculations above give
  $\pi^*\Phi\cdot|J\pi|\gtrsim|z|^{\delta}|w|^{\delta'}$
  where $\delta,\delta'>-1$. 
  Thus $(\pi^*\Phi\cdot|J\pi|)^2$ is locally integrable at $p$,
  completing the proof when $u$ has logarithmic singularities.

  We now consider the case of an arbitrary psh function. 
  As in the proof of Theorem~4.2 in~\cite{DK} 
  we use Demailly approximation to reduce to the preceding case.
  By Section~\ref{S1} 
  there exists a small neighborhood $B'$ of the origin, and,
  for each $n>0$, finitely many holomorphic functions 
  $\{g_{nl}\}_{l=1}^k$ with $\int_{B'}|g_{nl}|^2\exp(-2nu)\le 1$, 
  such that the psh function $u_n=(2n)^{-1}\log\sum|g_{nl}|^2$
  (which has logarithmic singularities) 
  satisfies $|\nu(u)-\nu(u_n)|\le n^{-1}A(\nu)$
  for any $\nu\in\cVqm$. 
  Thus H\"older's inequality with $p=n$ and $q=n/(n-1)$ gives
  \begin{multline*}
    \int_{B'}|\psi|^2\exp(-2u)
    =\int_{B'}\left(\sum_1^k|g_{nl}|^2\exp(-2pu)\right)^{1/p}  
    \left(|\psi|^{2q}(\sum_1^k|g_{nl}|^2)^{-q/p}\right)^{1/q}\\ 
    \le k^{1/p}\left(
      \int_{B'}|\psi|^{2q}(\sum_1^k|g_{nl}|^2)^{-q/p}
    \right)^{1/q} 
    =k^{1/p}\left(\int_{B'}|\psi|^{2q}\exp(-2qu_n)\right)^{1/q}.
  \end{multline*}
  For any quasimonomial valuation $\nu$, we have
  \begin{equation*}
    |\chi_{u,\psi}(\nu)-\chi_{qu_n,\psi^q}(\nu)|
    =\left|
      \frac{\nu(u)}
      {\nu(\psi)+ A(\nu)}
      -\frac{q\nu(u_n)}
      {q\nu(\psi)+A(\nu)}
    \right|
    \lesssim \frac{1}{n}.
  \end{equation*}
  By assumption $\sup\chi_{u,\psi}<1$, hence 
  $\sup\chi_{qu_n,\psi^q}<1$ for large $n$.
  By the preceding argument, $|\psi|^{2q}\exp(-2qu_n)$ is locally
  integrable, thus so is $|\psi|^2\exp(-2u)$.
\end{proof}
%
%
%
%
\section{Singularity exponents and Kiselman numbers}\label{sec-arnold}
We now apply our machinery to study the singularities of
psh functions. 

Theorem~\ref{T703} immediately implies
that the complex singularity exponent $c(u)$
and Arnold multiplicity $\lambda(u)$ of a psh function $u$
equal those of its tree transform $g_u$
as defined by~\eqref{e501}.
In view of Remark~\ref{R101} and Lemma~\ref{lem-attain} we then have
\begin{Cor}\label{cor1}
  For any psh function $u$ we have 
  \begin{equation}\label{e401}
    \lambda(u)=\sup_{\nu\in\cVqm}\frac{\nu(u)}{A(\nu)}
    \qand
    c(u)=\inf_{\nu\in\cVqm}\frac{A(\nu)}{\nu(u)},
  \end{equation}
  and the supremum and infimum are attained at a valuation 
  of multiplicity one.

  Further, if $c=c(u)$, then $\exp(-2cu)$ is not locally integrable
  at the origin. In fact, $\vol\{u<\log r\}\gtrsim r^{2c(u)}$
  as $r\to0$.
\end{Cor}
\begin{Remark}
  This corollary contains an affirmative answer to the
  openness conjecture for psh functions in dimension $2$: 
  see~\cite[Remarks~5.3 and~4.4]{DK}.
\end{Remark}
\begin{Example}
  If $\psi=y^n+x^m$ with $m<n$, $m,n$ relatively prime,
  and $u=\log|\psi|$ then the supremum in~\eqref{e401} 
  must be attained at a valuation of the form $\nu_{y,t}$, 
  $1\le t\le\infty$.
  Thus $\lambda(\log|\psi|)=\sup_t\max\{nt,m\}/(1+t)=nm/(n+m)$.
\end{Example}
Let us rephrase Corollary~\ref{cor1} in terms of Kiselman numbers,
using Section~\ref{S10}. If $m(\nu)=1$, then 
$\nu=\nu_{y,t}$ where $(x,y)$ are local coordinates and $t\ge1$. 
Moreover, $\nu(u)=\nu^{x,y}_{t,1}(u)$ and $A(\nu)=1+t$.
Hence $\nu(u)/A(\nu)=\nu^{x,y}_{t/(1+t),1/(1+t)}(u)$.
\begin{Cor}\label{cor3}
  The Arnold multiplicity of $u$ is the supremum of all Kiselman 
  numbers $\nu^{x,y}_{a,b}(u)$ over all choices of local coordinates 
  $(x,y)$ and all choices of weights $(a,b)$ with $a+b=1$.
\end{Cor}
By fixing the local coordinates $(x,y)$ and letting the
weights vary we obtain a lower bound for the Arnold multiplicity
proved by Kiselman~\cite{kis2} (in any dimension).
We can also bound it in terms of the Lelong number $\nuL(u)$:
\begin{Cor}\label{C302}
  The Arnold multiplicity $\lambda(u)$ of $u$ satisfies 
  \begin{itemize}
  \item[(i)]
      $\frac12\nuL(u)\le\lambda(u)\le\nuL(u)$;
    \item[(ii)]
      $\lambda(u)=\nuL(u)$ iff $dd^cu$ 
      is the current of integration on a 
      smooth curve plus a current with zero Lelong number;
    \item[(iii)]
      if $\lambda(u)=\frac12\nuL(u)$, then
      the Lelong number of
      the strict transform of $u$ under a single blowup $\pi$ of the 
      origin is at most $\frac12\nuL(u)$ at any point on the
      exceptional divisor $\pi^{-1}(0)$.
  \end{itemize}
\end{Cor}
\begin{Remark}
  The bounds in~(i) are due (in any dimension) 
  to Skoda~\cite{skoda}. The characterization 
  in~(ii) sharpens a recent result by
  Blel and Mimouni~\cite{BM} (see also~\cite{mim}),
  who proved that if $u$ has Lelong number one, then
  $\exp(-2u)$ is locally integrable unless $dd^cu$ puts mass
  on an analytic curve.
  The implication in~(iii) seems to be new.
\end{Remark}
\begin{Remark}
  The converse to~(iii) is false, as is shown by the example
  $u=\log\max\{|y|,|x|^{1+\e}\}$, where $0<\e<1$. 
  Here $\nuL(u)=1$, $\lambda(u)=(1+\e)/(2+\e)>1/2$, but the 
  strict transform of $u$ has Lelong number zero at all points
  of $\pi^{-1}(0)$, but one, at which the Lelong number is 
  $\e$.
\end{Remark}
\begin{proof}[Proof of Corollary~\ref{C302}]
  All of this is a consequence of 
  Theorem~\ref{T703} and Proposition~\ref{P710}. 
  Indeed~(i) follows immediately, as does~(iii) in view 
  of Section~\ref{S11}.
  As for~(ii) we get that the tree measure $\rho_u$ is a point mass
  at a curve valuation, which must be 
  associated to a smooth analytic curve $\{\phi=0\}$. 
  Siu's Theorem now yields $u=\log|\phi|+u'$ with 
  $u'$ psh. Clearly $u'$ has zero Lelong number.
\end{proof}
%
%
%
%
\section{The ascending chain condition.} \label{sec-ACC}
We now wish to give a new proof of a result describing
the structure of the set of complex singularity exponents 
of holomorphic functions. 

As our approach is algebraic, we work in the 
general case, that is, 
$R$ is an equicharacteristic zero, two-dimensional,
regular local ring with algebraically closed residue field.
If $\psi\in R$, let $c(\psi)=c(g_\psi)$ and $\lambda(\psi)=c(\psi)^{-1}$
be the singularity exponent and Arnold multiplicity of $\psi$.
Here $g_\psi$ is the tree potential of $\psi$. Thus
$\lambda(\psi)=\sup_{\cVqm}\chi_\psi(\nu)$, where 
$\chi_\psi(\nu)=\nu(\psi)/A(\nu)$.
In the analytic case $R=\cO_0$ 
we have $c(\psi):=c(\log|\psi|)$ and $\lambda(\psi)=\lambda(\log|\psi|)$.

To motivate the result, first consider 
$\psi\in\fm$ irreducible and let $(\nu_i)_0^g$ be the
approximating sequence of $\nu_\psi$. It
follows from Corollary~\ref{cor1} that
the supremum of $\chi_\psi$ is attained at $\nu_1$.
This gives $c(\psi)=\frac{1}{\nu_1(\psi)}+\frac{1}{m(\psi)}$.
By Section~\ref{S3}, $\nu_1(\psi)$ is an integer, 
so that $c(\psi)\in\frac{1}{\N^*}+\frac{1}{\N^*}$,
a result that seems to have first been proved by 
Igusa~\cite{Igu77}; see also~\cite{Kuw99}.
When $\psi$ is no longer assumed to be irreducible,
$c(\psi)\not\in\frac{1}{\N^*}+\frac{1}{\N^*}$ in general
(see~\cite{phong-sturm} for a counter-example). Nevertheless we have
\begin{Thm}{(\cite[1992]{sho},~\cite[1999]{Kuw99},~\cite[2000]{phong-sturm})}
\label{ACC}
  The set $\mathbf{c}=\{c(\psi)\ ;\ \psi\in\fm\}$ satisfies the
  ascending chain condition (ACC): any increasing sequence in
  $\mathbf{c}$ is eventually stationary.
  
  Further, the limit points of $\mathbf{c}$ are $0$ and the rational
  numbers $1/a$ for $a\ge 1$.
\end{Thm}
\begin{proof}
  We will show that
  $\boldsymbol{\lambda}:=\{\lambda(\psi)\ ;\ \psi\in\fm\}$ satisfies
  the \emph{decreasing} chain condition (DCC) and that its
  limit points are $\N\cup\{\infty\}$.
  Our proof mainly consists of a careful study of the functions
  $\chi_\psi$ on $\cVqm$, for $\psi\in\fm$.

  Fix $\psi$ and let $\psi=\prod\phi_k^{a_k}$ be the decomposition
  into irreducible factors. For any $k$ we define the valuation 
  $\nu_k=\max\{\nu\in\cV\ ;\ m(\nu)=1,\ \nu\le\nu_{\phi_k}\}$. 
  This is either $\nu_{\phi_k}$ or a divisorial valuation. 
  Write $\nu_{kl}=\nu_k\wedge\nu_l$ and $\a_{kl}=\a(\nu_{kl})$.
  
  \medskip
  \textbf{Fact 1}. 
  For all $k,l$, either $\a_{kl}m_k$ or $\a_{kl}m_l$ is an 
  integer where $m_k=m(\phi_k)$.

  \smallskip
  We first note that $m_k \a_k$ is always an integer 
  (if $\a_k<\infty$). For this we refer repeatedly to
  Section~\ref{S3}. 
  As $\nu_{\phi_k}>\nu_k$, the generic
  multiplicity $b_k$ of $\nu_k$ divides $m_k$.
  Thus $\a(\nu_k)b(\nu_k)=\a_kb_k\in\N$, 
  hence $\a_km_k\in\N$. This shows Fact~1 
  when $\nu_{kl}=\nu_k$ or $\nu_l$. 
  If this is not the case, 
  then $\nu_{kl}<\nu_k$, $\nu_{kl}<\nu_l$.  
  As $m(\nu_k)=m(\nu_l)=1$, $b(\nu_{kl})=1$ and
  again $\a_{kl}\in\N$.

  \medskip
  \textbf{Fact 2}. 
  The supremum of $\chi_\psi$ is attained at a valuation
  $\nu_{kl}$ for some $k,l$.  

  \smallskip
  After unwinding definitions, Lemma~\ref{lem-attain}~(i) and~(ii)
  imply that the supremum of $\chi_\psi$ is attained either at
  $\nu_\fm$ or at one of the $\nu_{kl}$.
  If $\nu_{kl}=\nu_\fm$ for some $k$ and  $l$, the proof is complete,
  so assume that $\nu_{kl}>\nu_\fm$ for all $k$, $l$.
  Then there exists a unique tangent vector
  $\vv$ at $\nu_\fm$ such that
  $\nu_{kl}\in U(\vv)$ for all $k$, $l$.
  It is then straightforward to verify
  that $D_\vv\chi_\psi>0$. 
  Hence the maximum is not attained at $\nu_\fm$ and we are done.

  \medskip
  \textbf{Fact 3}. 
  For any $k,l$, we have 
  $\nu_{kl}(\psi)\in\frac{\N^*}{m_k}\cup\frac{\N^*}{m_l}$.

  \smallskip
  If $\nu_{\phi_i}\ge\nu_{kl}$, we have $\nu_{kl}(\phi_i)=\a_{kl}m_i$ 
  which belongs to $\N^*/m_k \cup \N^*/m_l$ by Fact~1. 
  If $\nu_{\phi_i}\not\ge\nu_{kl}$, then 
  $\nu_{kl}(\phi_i)=\nu_{ik}(\phi_i)$. Then
  either $\nu_{ik}=\nu_i$, and $\nu_{kl}(\phi_i)\in\N^*$; 
  or $\nu_{ik}<\nu_i$, $\nu_{ik}<\nu_{kl}$ and $\a_{ik}\in\N^*$ so
  that $\nu_{kl}(\phi_i)\in\N^*$.
  
  \medskip
  \textbf{Fact 4}. 
  Suppose $\lambda(\psi)=\chi_\psi(\nu_{kl})$.  
  Set $I=\{i\ ;\ \nu_{\phi_i}\ge\nu_{kl}\}$ 
  and let $J$ be its complement. Then
  \begin{equation}\label{ineq}
    \sum_J a_j\nu_{kl}(\phi_j)\le\sum_Ia_im_i.
  \end{equation}
  This inequality follows from the assumption 
  $\lambda(\psi)=\chi_\psi(\nu_{kl})$.
  For any $t$ in an interval $T=(\a_{kl}-\eta,\a_{kl})$ for $\eta$
  sufficiently small, we have 
  \begin{equation*}
    \chi(\nu_{\psi_i,t})
    =\frac{t(\sum_Ia_im_i)+\sum_Ja_j\nu_{kl}(\phi_j)}{1+t}.
  \end{equation*}
  As the supremum of $\chi$ is attained when $t = \a_{kl}$, we infer
  that $\chi' \ge 0$ on $T$. A direct computation then shows~\eqref{ineq}.
  
  \medskip
  We are now able to prove the theorem. First note that it suffices
  prove that $\boldsymbol{\lambda}\cap(0,C)$ satisfies the DCC
  for any $C>0$. So pick a sequence $\psi_n\in\fm$ such that
  $\lambda(\psi_n)$ is decreasing and bounded from above by $C$. 
  Since $\lambda(\psi)=\sup\nu(\psi)/A(\nu)\ge m(\psi)/2$, 
  we can bound the multiplicity $m(\psi_n)\le 2C$. 
  For any $n$, we introduce $\phi_k^n$, $\a_{kl}^n$, $I^n$, $J^n$ as above. 
  To simplify notation we drop the superscript $n$. 
  We also denote by $\cN=\{\frac{p}{q}\ ;\ p,q\in\N^*,\ q\le2C\}$. 
  This is a discrete semi-group. Fact~2 implies that we can write
  \begin{align*}
    \lambda(\psi_n)
    &=\frac{(\sum_I a_i m_i)\a_{kl}+\sum_Ja_j\nu_{kl}(\phi_j)}{1+\a_{kl}}\\
    &=(\sum_Ia_i m_i)-\frac{1}{1+\a_{kl}}\left[
      \sum_Ia_im_i-\sum_Ja_j\nu_{kl}(\phi_j)
    \right].
  \end{align*}
  We have $\sum_Ia_im_i\le m(\psi_n)\le2C$, so we may
  assume that $\sum_Ia_im_i$ is constant, say equal to $D$,
  for all $n$. Fact~3 implies that $\sum_Ja_j\nu_{kl}(\phi_j)$
  belongs to the discrete set $\cN$. It is also bounded by $D$ 
  by Fact~4, so we may assume it is constant.
  Thus $\lambda(\psi_n)=D-\frac{E}{1+\a_{kl}}$ for some $E$.
  Now $\lambda(\psi_n)$ is increasing, hence $\a_{kl}$ is decreasing. 
  Finally Fact~1 shows that $\a_{kl}\in\cN$ and therefore
  is constant for large $n$. This proves that
  $\boldsymbol{\lambda}$ satisfies the DCC.
  
  In order to find the limit points of $\boldsymbol{\lambda}$,
  suppose $\lambda(\psi_n)\to \lambda\in\boldsymbol{\lambda}$, 
  $\lambda\ne\infty$. 
  By extracting a subsequence and using the DCC 
  we may assume that $\lambda(\psi_n)$ is increasing.
  The argument above applies again: 
  for $n$ large, we have $\lambda(\psi_n)=D-E(1+\a_{kl})^{-1}$ 
  for some constants $D,E>0$, where $D$ is a integer. 
  If $\lambda(\psi_n)$ is not stationary, $\a_{kl}$ increases to infinity,
  so that $\lambda=D\in\N$. 
  This shows that the limit points of $\boldsymbol{\lambda}$ are included 
  in $\N\cup\{\infty\}$. 
  Conversely, we have
  $\lambda(y^m+x^n)=mn/(m+n)$. Letting $m\to\infty$ while keeping
  $n$ constant gives $n\in\boldsymbol{\lambda}$.
  Letting $m,n\to\infty$ gives $\infty\in\boldsymbol{\lambda}$.
  This concludes the proof.
\end{proof} 
%
%
%
%
\section{Ideals as multiplier ideals}\label{sec-asmult}
Proposition~\ref{P708} asserts that the multiplier ideal associated to
any tree potential is integrally closed.
Our aim is to prove the converse statement, thereby
giving a new proof of a theorem by Lipman and Watanabe~\cite{LipWat03}.
\begin{Thm}\label{T101}
  If $J$ is an integrally closed ideal of $R$, then
  $J=\cJ(g)$ for some tree potential $g$.
\end{Thm}
\begin{Remark}\label{R704}
  Our proof is constructive and generates a tree potential of the form
  $g=c\,g_{I}$, where $I$ is an integrally closed ideal and $c>0$. 
  The choices of $I$ and $c$ are not unique.
\end{Remark}
\begin{Remark}\label{R805}
  Lipman and Watanabe work on more general rings.
  They also prove that if $J=I_\nu$ with $\nu$ divisorial, then 
  there exists an \emph{ideal} $I$ with $\cJ(I)=J$ iff $b(\nu)=1$.
  More generally it would be interesting to know what integrally
  closed ideals are of the form $\cJ(I)$ for an ideal $I$.
\end{Remark}
\begin{proof}
  Any integrally closed ideal in $R$ 
  is the product of a principal ideal and a primary ideal.
  In view of Proposition~\ref{P304} we may hence assume
  that $J$ is primary.
  
  Let us recall some notation from the analysis of tree transforms
  of ideals in Section~\ref{S5}.
  The tree measure of $J$ on $\cV$ is of the form
  $\rho_J=\sum_1^rn_ib_i\nu_i$.
  Here $n_i\in\N$ and $\nu_i$ is divisorial with generic 
  multiplicity $b_i$.
  Write $\a_i=\a(\nu_i)$ and $A_i=A(\nu_i)$.
  Set $m=m(J)=\sum_in_ib_i$ and
  let $\cT_J$ be the support of the tree potential $g_J$, 
  \ie the smallest subtree containing $\nu_\fm$ and all
  the $\nu_i$.
  
  \begin{figure}[ht]
    \includegraphics[width=0.9\textwidth]{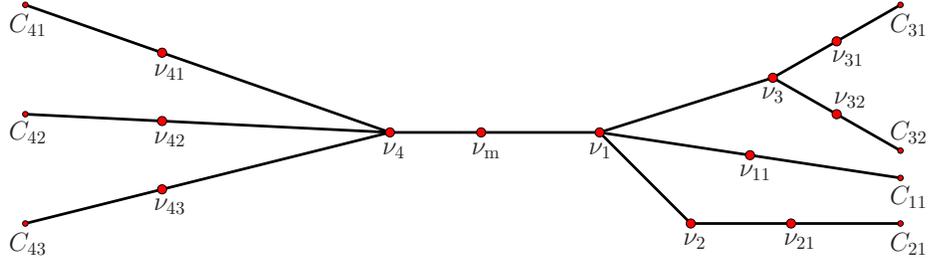}
    \caption{Construction of a tree potential with a predetermined
      multiplier ideal $I$. The Zariski factorization of $I$
      is $I=I_1I_2I_3^2I_4^3$ and $\nu_i$ is the Rees valuation
      of $I_i$, $1\le i\le 4$. The curves $C_{ij}$ are chosen
      as curvettes of $\nu_i$ and the valuations $\nu_{ij}$
      are well-chosen on the segments $[\nu_i,\nu_{C_{ij}}]$.
      See the proof of Theorem~\ref{T101}.}\label{F3}
  \end{figure}

  Pick $\psi_{ij}\in\fm$ irreducible for $1\le i\le r$, $1\le j\le n_i$ 
  such that $m(\psi_{ij})=b_i$, $\nu_{\psi_{ij}}>\nu_i$ and 
  the $\nu_{\psi_{ij}}$ represent $n_i$ distinct tangent vectors at 
  $\nu_i$. The curves $C_{ij}:=\{\psi_{ij}=0\}$ are then curvettes 
  of $\nu_i$.
  Define $\nu_{ij}$ to be the unique
  quasimonomial valuation in $]\nu_i,\nu_{\psi_{ij}}[$ 
  whose thinness
  $A_{ij}=A(\nu_{ij})$ is given by $A_{ij}-A_i=R_i$ for a constant
  $R_i>0$ depending only on $i$ to be fixed shortly.
  See Figure~\ref{F3}.
  Write $\a_{ij}=\a(\nu_{ij})$ so that $b_i(\a_{ij}-\a_i)=R_i$. 
  For a given $\e>0$, we define the tree potential $g$ by
  \begin{equation*} 
    \Delta g=\left(1+\e\right)\sum_{i\le r,\ j\le n_i}b_i\nu_{ij}.  
  \end{equation*}
  Let $\cT_g$ be the support of $g$ (\ie $\cT_J$ with the 
  segments $[\nu_i,\nu_{ij}]$ attached)
  and let $c(J)=\inf_\nu A(\nu)/\nu(J)$ be the complex 
  singularity exponent of the ideal $J$ (see Corollary~\ref{cor1}). 
  In the sequel, we fix $\e<m(J)^{-1}$. 
  With this choice, we now show
  that we can pick $R_i>0$, for $i\le r$, so that
  \begin{equation}\label{E101} 
    \nu_i(J)+A_{ij}-b_i^{-1}
    <g(\nu_{ij})
    <\nu_i(J)+A_{ij}.
  \end{equation} 
  To see this, notice that
  \begin{equation}\label{E100} 
    g(\nu_{ij})= (1+\e)(\nu_i(J)+R_i).  
  \end{equation} 
  Hence $g(\nu_{ij})-\nu_i(J)-A_{ij}=\e\nu_i(J)-A_i+\e R_i=:S_i$.
  We want to pick $R_i>0$ so that $S_i\in(-b_i^{-1},0)$.
  Clearly $S_i$ is an affine function of $R_i$ and 
  $S_i\to\infty$ as $R_i\to\infty$. Hence it suffices to show that
  $S_i<0$ when $R_i=0$ \ie that $\e\nu_i(J)<A_i$.
  But the tree potential properties of $g_J$ give
  $\nu_i(J)\le m\a_i$. Hence $\e\nu_i(J)<\a_i<A_i$.

  Finally, we note that $g$ is proportional to the tree transform of
  some integrally closed ideal. Indeed multiplying 
  $\sum b_i\nu_{ij}$
  by the least common multiple of all the generic multiplicities
  $b(\nu_{ij})$ we obtain the tree measure $\rho_I$ of some
  ideal $I$.
  We can thus write $g=c\,g_I$ for some $c>0$ and some
  integrally closed ideal $I$.

  \smallskip
  We claim that $\cJ(g)=J$. 
  Let us first prove that $\cJ(g)\supset J$.
  By Proposition~\ref{P301}~(ii) it suffices to show that 
  $\nu(J)+A(\nu)>g(\nu)$ for $\nu\in\cT_J$ and 
  $\nu=\nu_{ij}$.
  The latter case is taken care of by~\eqref{E301}.
  As for the former, if $\nu\in\cT_J$, then $g(\nu)=(1+\e)\nu(J)$.
  Moreover, $\nu(J)\le m(J)\a(\nu)<m(J)A(\nu)$, so
  \begin{equation*}
    \frac{g(\nu)}{\nu(J)+A(\nu)}
    =\frac{1+\e}{1+A(\nu)/\nu(J)}
    <\frac{1+\e}{1+m(J)^{-1}}<1.
  \end{equation*}

  We now complete the proof by proving that $\cJ(g)\subset J$.
  Fix $\psi\in\cJ(g)$. Then  
  \begin{equation}\label{E200}
    \nu(\phi)>g(\nu)-A(\nu)
  \end{equation} 
  for all $\nu\in\cVqm$. 
  By Section~\ref{S5}
  we must show that
  $\nu_i(\phi)\ge \nu_i(J)$ for all $i\le r$.

  For $i\le r$ and any $j\le n_i$,~\eqref{E101}
  and~\eqref{E200} imply that
  \begin{equation*}
    \nu_{ij}(\phi)
    >g(\nu_{ij})-A_{ij}
    >\nu_i(J)-b_i^{-1}.
  \end{equation*}
  We let $\cN\subset\{1,\dots,r\}$ be the set of
  $i$ for which there exists $j$ such that
  $\phi$ has no irreducible factor $\tilde{\phi}$ with
  $\nu_{\tilde{\phi}}\wedge\nu_{ij}>\nu_i$.
  After relabeling we have $\cN=\{1,\dots,k\}$.
  
  For $1\le i\le k$, we infer immediately that
  $\nu_i(\phi)=\nu_{ij}(\phi)>\nu_i(J)-b_i^{-1}$.
  Since $\nu_i$ is divisorial we have
  $\nu_i(\phi)\in b_i^{-1}\N$ (see Section~\ref{S6}),
  hence $\nu_i(\phi)\ge\nu_i(J)$.
  
  If $k=r$, the proof is complete. 
  Thus assume $k<r$ and pick $i$ with $k<i\le r$. 
  For any $j\le n_i$, $\phi$ admits an
  irreducible factor $\phi_{ij}$ with 
  $\nu_{\phi_{ij}}\wedge\nu_{ij}>\nu_i$. 
  The multiplicity of $\phi_{ij}$ is necessarily a 
  multiple of $b_i$, say $p_{ij}b_i$. 

  We now replace $\phi_{ij}$ by $\psi_{ij}^{p_{ij}}$, where 
  $\psi_{ij}$ are defined above. By doing this replacement 
  whenever $k<i\le r$ and $1\le j\le n_i$ we obtain 
  from $\phi$ a new $\hphi\in\fm$, 
  for which $\nu(\hat{\phi})=\nu(\phi)$ for all $\nu\in\cT_J$.

  We can then write
  $\hphi=\phi'\,\prod_{k<i\le r}\psi_i$
  where $\psi_i=\prod_j\psi_{ij}$ and $\phi'\in\fm$.
  Notice that if $\nu\in\cT_J$ then
  $\nu(\psi_i)=b_i\,\nu\cdot\nu_i$ and 
  $\nu(J)=\sum_1^r n_ib_i\,\nu\cdot\nu_i$.
  Thus, for $\nu\in\cT_J$ we have
  \begin{equation}\label{e708}
    \nu(\phi)-\nu(J)
    =\nu(\hat{\phi})-\nu(J)
    =\nu(\phi')-\sum_1^kn_ib_i\, \nu\cdot\nu_i
    =\nu(\phi')-\nu(J'),
  \end{equation}
  where $J'$ is the integrally closed ideal associated to
  the measure $\sum_1^kn_ib_i\nu_i$ (see Section~\ref{S5}).
  We have seen that 
  $\nu_i(\phi')-\nu_i(J')=\nu_i(\phi)-\nu_i(J)>-b_i^{-1}$ 
  for all $i\le k$. As above we conclude that in fact
  $\nu_i(\phi')-\nu_i(J')\ge0$. As the support of $\rho_{J'}$ is included
  in the $\nu_i$'s, we infer
  $\nu(\phi')\ge\nu(J')$ for all $\nu\in\cVqm$. 
  Thus~\eqref{e708} gives $\nu(\phi)\ge\nu(J)$ 
  for all $\nu\in\cT_J$. 
  Hence $\phi\in J$, which concludes the proof.
\end{proof}
\begin{Remark}
  It follows from the proof that we may in fact choose the
  numbers $R_i$ of the form $R_i=n/b_i$ for $n\in\N^*$. 
  One can then check that our construction coincides with that
  of Lipman and Watanabe~\cite{LipWat03}.
\end{Remark}
%
%
%
%
\section{Equisingularity}\label{sec-equi}
In general we may refer to two objects
as \emph{equisingular} if they have the same tree transform.
Two ideals are equisingular iff
they have the same integral closure.
(This is classical, but see~\cite[Theorem~8.12]{valtree}.)
Two psh functions are equisingular iff
their pullbacks by any composition of blowups have the same
Lelong number at all points on the exceptional divisor,
see~\cite[Proposition~6.2]{pshsing}.

Our objective now is to show that equisingularity can be detected 
at the level of multiplier ideals.
This is achieved through the following theorem, which allows us
to recover a tree potential from the 
multiplier ideals of all its constant multiples.
\begin{Thm}\label{thm-apx-gen}
  Suppose $h$ is a potential. For $t>0$, let 
  $h_t$ be the tree transform of the ideal $\cJ(th)$.
  Then $h-t^{-1}A\le t^{-1}h_t\le h$, where $A$ denotes
  thinness. In particular $t^{-1}h_t$ converges to $h$
  in $\cP$ as $t\to\infty$.
\end{Thm}
\begin{Remark}
  The sequence $(n^{-1}h_n)_{n=1}^\infty$ 
  can be thought of as an analogue for a tree potential
  of the Demailly approximating sequence of a psh function.
\end{Remark}
\begin{Cor}\label{cor-apx1}
  Let $h_1,h_2\in\cP$. Then $h_1=h_2$
  iff $\cJ(th_1)=\cJ(th_2)$ for all $t\ge 0$.
\end{Cor}
Using the results of Sections~\ref{sec-ideal} and~\ref{sec-psh}
we immediately infer two corollaries.
\begin{Cor}\label{cor-apx3}
  Let $I_1,I_2$ be integrally closed ideals.
  Then $I_1=I_2$ iff $\cJ(I_1^n)=\cJ(I_2^n)$ for all $n\in\N^*$.
\end{Cor}
\begin{Cor}\label{cor-apx2}
  Let $u_1,u_2$ be psh functions. 
  Then $\mathcal{J}(t\,u_1)=\mathcal{J}(t\,u_2)$ for all $t\ge 0$, 
  iff the tree transforms of $u_1$ and $u_2$ coincide.
\end{Cor}
\begin{Remark}
  Corollary~\ref{cor-apx2} can be deduced directly from
  Ohsawa-Takegoshi's theorem. Similarly, Corollary~\ref{cor-apx3} is a
  consequence of Skoda's theorem, which implies 
  $\cJ(I^n)=I^{n-1}\cJ(I)$.
\end{Remark}
\begin{Remark}
  It is not true that if $I_\bullet$, $J_\bullet$ are graded systems 
  of ideals with $\cJ(I_\bullet^c)=\cJ(J_\bullet^c)$ for all $c>0$,
  then $I_\bullet=J_\bullet$. See Example~\ref{E502}.
\end{Remark}
We give a proof of Theorem~\ref{thm-apx-gen}
based on the subadditivity property for multiplier ideals proved
in~\cite{DEL}.
\begin{proof}[Proof of Theorem~\ref{thm-apx-gen}]
  One inequality is elementary:
  for any $\phi\in\cJ(th)$ we have $\nu(\phi)>t\,h(\nu)-A(\nu)$,
  hence $t^{-1}h_t\ge h-A/t$.
  The other inequality, on the other hand, is quite deep, but it 
  is an immediate consequence of the following result,
  whose proof is given below.
\end{proof}
\begin{Lemma}\label{osh-tak}
  For any tree potential $h$, the tree transform of $\cJ(h)$ 
  is dominated by $h$. In other words, $g_{\cJ(h)}\le h$.
\end{Lemma}
\begin{Remark}
  Explicitly, this lemma asserts that
  for any tree potential $h$, and any $\nu\in\cVqm$, one can find 
  $\phi\in\cJ(h)$ so that $\nu(\phi)\le h(\nu)$. 
  It would be interesting to have a more direct construction
  of such an element $\phi$.
\end{Remark}
\begin{proof}[Proof of Lemma~\ref{osh-tak}]
  Let $\rho=\rho_h$ be the measure associated to $h$. We will
  gradually build the proof using successively more general measures 
  $\rho$.
  
  First suppose $h=g_I$ for an integrally closed ideal $I$.
  If  $\psi\in I$, then $\nu(\psi)\ge\nu(I)$ for all $\nu\in\cV$. 
  Hence $\nu(\psi)+A(\nu)>\nu(I)$ for all $\nu\in\cVqm$
  and $\rho_\psi\{\nu\}+m(\nu)>\nu_I\{\nu\}$ for all
  curve valuations $\nu$. 
  By Proposition~\ref{P301} this implies $\psi\in\cJ(I)$.
  Hence $I\subset \cJ(I)$, so Lemma~\ref{osh-tak} is
  immediate in this case.
  
  We now consider the more general case $h=c\,g_I$, 
  where $I$ is an integrally closed ideal and $c>0$ is rational.
  By the subadditivity theorem of
  Demailly-Ein-Lazarsfeld~\cite{DEL} we have
  $\cJ(I^{cq})\subset\cJ(I^c)^q$ for any $q$. 
  (Our situation is strictly speaking not contained in~\cite{DEL};
  a more general statement than what we need can be found 
  in~\cite{TakWat03}.)
  Therefore $I^{cq}\subset\cJ(I^c)^q$ if $cq\in\N$.
  Hence $cq\,g_I\ge q\,g_{\cJ(I^c)}$, which gives $h\ge g_{\cJ(h)}$.
  
  The next more general case is when the
  support of $\rho$ is included in a finite tree
  $\cT$ whose ends are quasimonomial valuations. 
  We approximate $\rho$ by atomic measures
  $\rho_n=\sum_i c^n_i\nu^n_i$. 
  We can suppose $\nu_i^n$ are all divisorial valuations, and $c^n_i$ 
  are rational numbers. We can also impose that $h_n=g_{\rho_n}$ decreases
  towards $g_\rho=h$. To do so, note that by linearity, we may take the
  support of $\rho$ to be included in a segment $T=[\nu_\fm,\nu_0]$.
  The set of divisorial valuations in $T$ coincides with the set of
  valuations with rational skewness, and is hence dense in $T$. 
  We now have to approximate a concave function on $T$ by a decreasing
  sequence of piecewise linear functions with rational slopes, 
  and whose second derivatives are measures supported on the rationals. 
  This can be done in an elementary way.
  
  By Section~\ref{S5}, $q\,h_n$ is the tree transform of an integrally
  closed ideal $I$, if $q$ is the product of the denominators of $c^n_i$
  with the generic multiplicities of $\nu^n_i$ over all $i$.  By the
  preceding argument, the tree transform of $\cJ(h_n)$ is bounded
  by $h_n$. As $h\le h_n$ we infer $g_{\cJ(h)}\le g_{\cJ(h_n)}\le h_n$. 
  We conclude, by letting $n\to\infty$, that $g_{\cJ(h)}\le h$ as desired.
  
  Finally consider a general tree potential $h$.
  First assume that $\rho\{\nu\}<m(\nu)$ for all curve valuations 
  $\nu$. Since $\rho$ has finite mass we may pick
  $\e>0$ such that in fact $\rho\{\nu\}<(1-\e)m(\nu)$.
  Define $\cT=\{\nu\in\cV\ ;\ \rho\{\mu\ge\nu\}\ge (1-\e)m(\nu)\}$.
  By Lemma~\ref{L711}, $\cT$ is a finite subtree of 
  $\cV$ and $\cJ(h_\cT)=\cJ(h)$.
  Our choice of $\e$ means that the ends of $\cT$ 
  are quasimonomial. Hence the previous analysis gives
  $g_{\cJ(h_\cT)}\le h_\cT$. Since $h_\cT\le h$ this yields
  $g_{\cJ(h)}\le h$.

  In the most general case we can decompose 
  $h=g_{\phi}+h'$, where
  $\phi\in R$ (possibly reducible) and
  $h'$ is a tree potential with associated measure
  satisfying $\rho'\{\nu\}<m(\nu)$ for all curve valuations
  $\nu$. By repeated use of Proposition~\ref{P304} we have
  $\cJ(h)=\phi\cJ(h')$.
  Therefore the previous step gives
  \begin{equation*}
    g_{\cJ(h)}
    =g_\phi+g_{\cJ(h')}
    \le g_\phi+h'
    =h,
  \end{equation*}
  which concludes the proof.
\end{proof}
%
%
%
%
\appendix
\section{Subadditivity and Skoda's Theorem}\label{sec-subadd}
In this appendix we prove two fundamental properties of
multiplier ideals associated to tree potentials. 
The proofs appear here as they rely on the corresponding 
properties of multiplier ideals of formal powers of ideals.
It would be interesting to have direct proofs, not
using a reduction to the case of powers of ideals.
\begin{proof}[Proof of Proposition~\ref{P302}]
  Write $h=h_1+h_2$. 
  Let $\rho_i=\rho_{h_i}$ and $\rho=\rho_{h}=\rho_1+\rho_2$
  be the measures associated to $h_i$ and $h$, respectively.
  As in the proof of Lemma~\ref{osh-tak}
  we will gradually build the proof using successively more 
  general measures.
  
  First suppose that $h_i=c_ig_{I_i}$, $i=1,2$,
  for ideals $I_1$, $I_2$ and $c_1,c_2>0$ rational. 
  By Proposition~\ref{P401} our subadditivity statement 
  translates into 
  $\cJ(I_1^{c_1}\cdot I_2^{c_2})\subset\cJ(I_1^{c_1})\cdot\cJ(I_2^{c_2})$
  which holds by~\cite{DEL} (see the proof of Lemma~\ref{osh-tak}).

  Next suppose the 
  support of $\rho$ is included in a finite tree
  $\cT$ whose ends are quasimonomial valuations. 
  As in the proof of Lemma~\ref{osh-tak} we 
  approximate $h_i$ from above by 
  tree potentials of the form $h_i^n=c_i^ng_{I_i^n}$, where
  $c_i^n>0$ are rational and $I_i^n$ are primary ideals.
  By semicontinuity (Proposition~\ref{P709}) we have 
  $\cJ(h_i^n)=\cJ(h_i)$ and
  $\cJ(h_1^n+h_2^n)=\cJ(h)$ for large $n$, completing the
  proof in this case.

  Now consider general tree potentials $h_i$. First assume
  that $\rho_i\{\nu\}<m(\nu)$ for all curve valuations 
  $\nu$, $i=1,2$. Since $\rho_i$ has finite mass we may pick
  $\e>0$ so that in fact $\rho_i\{\nu\}<(1-\e)m(\nu)$.
  Define $\cT_i=\{\nu\in\cV\ ;\ \rho_i\{\mu\ge\nu\}\}$.
  By Lemma~\ref{L711}, $\cT_i$ is a finite subtree of 
  $\cV$ and $\cJ((h_i)_{\cT_i})=\cJ(h_i)$.
  Our choice of $\e$ means that the ends of $\cT_i$ 
  are quasimonomial. Hence the previous analysis gives
  \begin{equation*}
    \cJ(h_1+h_2)
    \subset\cJ((h_1)_{\cT_1}+(h_2)_{\cT_2})
    \subset\cJ((h_1)_{\cT_1})\cdot\cJ((h_2)_{\cT_2})
    =\cJ(h_1)\cdot\cJ(h_2).
  \end{equation*}

  In the most general case we can decompose 
  $h_i=g_{\phi_i}+h'_i$, where
  $\phi_i\in R$ (possibly reducible) and
  $h'_i$ is a tree potential with associated measure
  satisfying $\rho'_i\{\nu\}<m(\nu)$ for all curve valuations
  $\nu$. By repeated use of Proposition~\ref{P304} we have
  $\cJ(h_i)=\phi_i\cJ(h'_i)$ and $\cJ(h)=\phi_1\phi_2\cJ(h'_1+h'_2)$. 
  Therefore the previous step gives
  \begin{equation*}
    \cJ(h_1+h_2)
    =\phi_1\phi_2\cJ(h'_1+h'_2)
    \subset\phi_1\phi_2\cJ(h'_1)\cdot\cJ(h'_2)
    =\cJ(h_1)\cdot\cJ(h_2),
  \end{equation*}
  which completes the proof.
\end{proof}
\begin{proof}[Proof of Proposition~\ref{P303}]
  The proof is based on a reduction to 
  the case $h=c\,g_I$ for an integrally closed ideal $I$
  and $c>0$ rational. In the latter case 
  the result is well-known (see~\cite[Theorem~9.6.21]{laz-book}) 
  in view of Proposition~\ref{P401}. 
  The reduction goes along exactly
  the same steps as the proof of subadditivity above.
  The details are left to the reader.
\end{proof}
%
%
%
%

\end{document}